\documentclass[12pt]{amsart}

\usepackage{epsfig, fullpage}
\usepackage[mathscr]{euscript}
\allowdisplaybreaks

\newtheorem{theorem}{Theorem}
\newtheorem{proposition}[theorem]{Proposition}
\newtheorem{lemma}[theorem]{Lemma}
\newtheorem{corollary}[theorem]{Corollary}

\newtheorem{definition}[theorem]{Definition}

\theoremstyle{definition}
\newtheorem{remark}[theorem]{Remark}
\newtheorem{example}[theorem]{Example}

 \newcounter{thlistctr}
 \newenvironment{thlist}{\
 \begin{list}%
 {\alph{thlistctr}}%
 {\setlength{\labelwidth}{2ex}%
 \setlength{\labelsep}{1ex}%
 \setlength{\leftmargin}{6ex}%
 \usecounter{thlistctr}}}%
 {\end{list}}

\newcommand{\ole}{{\,\sqsubseteq\,}}
\newcommand{\iole}{{\,\sqsubseteq_{\mathbf I}\,}}
\newcommand{\Excl}{{\sf Excl}}
\newcommand{\Primes}{{\sf Primes}}
\newcommand{\Even}{{\sf Even}}
\newcommand{\Odd}{{\sf Odd}}
\newcommand{\idom}{{\mathbb{IDOM}}}
\newcommand{\dom}{{\mathbb{DOM}}}
\newcommand{\Seq}{{\sf Seq}}
\newcommand{\wSeq}{\widehat{\sf Seq}}
\newcommand{\Set}{{\sf Set}}

\newcommand{\MSet}{{\sf MSet}}
\newcommand{\wMSet}{\widehat{\sf MSet}}
\newcommand{\Cyc}{{\sf Cycle}}
\newcommand{\DCycle}{{\sf DCycle}}
\newcommand{\wDCycle}{\widehat{\sf DCycle}}
\newcommand{\wCycle}{\widehat{\sf Cycle}}
\newcommand{\Cycle}{{\sf Cycle}}
\newcommand{\DSing}{{\sf DomSing}}
\newcommand{\Spec}{{\sf Spec}}

\newcommand{\sC}{{\sf C}}
\newcommand{\sD}{{\sf D}}
\newcommand{\sH}{{\sf H}}
\newcommand{\sS}{{\sf S}}

\newcommand{\cA}{{\mathcal A}}
\newcommand{\cC}{{\mathcal C}}
\newcommand{\cO}{{\mathcal O}}
\newcommand{\cT}{{\mathcal T}}
\newcommand{\cU}{{\mathcal U}}

\newcommand{\bA}{{\bf A}}
\newcommand{\bB}{{\bf B}}
\newcommand{\bD}{{\bf D}}
\newcommand{\bE}{{\bf E}}
\newcommand{\bF}{{\bf F}}
\newcommand{\bG}{{\bf G}}
\newcommand{\bH}{{\bf H}}
\newcommand{\bL}{{\bf L}}
\newcommand{\bP}{{\bf P}}
\newcommand{\bQ}{{\bf Q}}
\newcommand{\bR}{{\bf R}}
\newcommand{\bS}{{\bf S}}
\newcommand{\bT}{{\bf T}}
\newcommand{\bU}{{\bf U}}
\newcommand{\bV}{{\bf V}}
\newcommand{\bZ}{{\bf Z}}

\newcommand{\bbC}{{\mathbb C}}
\newcommand{\bbD}{{\mathbb D}}
\newcommand{\bbM}{{\mathbb M}}
\newcommand{\bbN}{{\mathbb N}}
\newcommand{\bbP}{{\mathbb P}}
\newcommand{\bbR}{{\mathbb R}}
\newcommand{\bbU}{{\mathbb U}}

\newcommand{\ro}{{\rm o}}

\title{Counting Rooted Trees\,:\\ The Universal Law
$t(n)\,\sim\,C \rho^{-n} n^{-3/2}$}

\author{Jason P. Bell}
\thanks{We are greatly indebted to the referee for bringing up 
  important questions, especially regarding the role of $\Set$, that led 
  us to thoroughly rework the paper.}
\address{Department of Mathematics, Simon Fraser University,8888
  University Dr., Burnaby, BC, V5A 1S6, Canada}
\email{jpb@math.sfu.ca}

\author{Stanley N. Burris}
\thanks{The second and third authors 
    would like to thank NSERC for support of this research.}
\address{Department of Pure Mathematics, University of Waterloo,Waterloo, Ontario, N2L 3G1, Canada}
\email{snburris@thoralf.uwaterloo.ca}
\urladdr{www.thoralf.uwaterloo.ca}

\author{Karen A. Yeats}
\address{Department of Mathematics and Statistics, Boston University,
  111 Cummington St, Boston, MA 02215, USA}
\email{kayeats@math.bu.edu}

\date{July 27, 2006}

\begin{document}

\begin{abstract}
Combinatorial classes $\cT$ that are recursively defined using 
combinations of
the standard  \emph{multiset}, \emph{sequence}, \emph{directed cycle} and \emph{cycle} 
constructions, and their restrictions, 
have generating series $\bT(z)$ with a positive radius of
convergence; for most of these a simple test can be used 
to quickly show that the form of the asymptotics is the same 
as that for the class of rooted trees: $C \rho^{-n} n^{-3/2}$\,,
where $\rho$ is the radius of convergence of $\bT$.
\end{abstract}

\maketitle

\section{Introduction}

The class of rooted trees, perhaps with additional structure as in the planar case,
is unique among the well studied classes of structures. It is 
so easy to find endless possibilities for defining interesting
subclasses as {\em the fixpoint of a class construction}, where the constructions
used are combinations of a few standard constructions like {\em sequence}, 
{\em multiset} and {\em add-a-root}.  This fortunate situation is based on a simple 
reconstruction property: removing the root from a tree gives a collection of 
trees (called a forest); and it is trivial to  reconstruct the 
original tree from the forest (by adding a root). 

Since we will be frequently referring to {\em rooted} trees, and rarely to 
{\em free} (i.e., {\em unrooted}) trees, from now on we will assume, unless 
the context says otherwise, that the word `tree' means `rooted tree'.

\subsection{Cayley's fundamental equation for trees}

Cayley \cite{Cy1} initiated the tree investigations\footnote{This was in 
the context of an algorithm for expanding partial differential operators.
Trees play an important role in the modern theory of differential equations
and integration---see for example Butcher \cite{bu1972}.}
in 1857 
when he presented the well known infinite product 
representation\footnote{This representation
uses $t(n)$ to count the number of trees on $n$ vertices. Cayley actually
used $t(n)$ to count the number of trees with $n$ edges, so his formula was
\[
\bT(z)\ =\ z \prod_{j\ge 1}\big(1-z^j\big)^{-t(j-1)}.
\]
}
\[
\bT(z)\ =\ z \prod_{j\ge 1}\big(1-z^j\big)^{-t(j)}\,.
\]
Cayley used this to calculate $t(n)$ for $1\le n\le 13$\,.  
More than a decade later (\cite{Cy3}, \cite{Cy4}, \cite{Cy6}) 
he used this method to give
\emph{recursion procedures} for finding the coefficients of generating 
functions for the chemical diagrams of certain families of compounds.  

\subsection{P\'olya's analysis of the generating series for trees}

Following on Cayley's work and further contributions by chemists, 
P\'olya published his classic 1937 paper\footnote{Republished in book form 
in \cite{P:R}.} 
that presents:
(1) his {\em group-theoretic} approach to enumeration, and (2) the primary analytic 
technique to establish the \emph{asymptotics} of recursively defined classes 
of trees. Let us review the latter as it has provided the paradigm for all 
subsequent investigations into generating series defined by 
recursion equations.

Let $\bT(z)$ be the generating series for the class of all unlabelled 
trees. P\'olya first converts Cayley's equation to the form
\[
\bT(z)\ =\ z\cdot\exp\Big(\sum_{m\ge 1}\bT(z^m)/m\Big).
\]
From this he quickly deduces that: the radius of convergence $\rho$
of $\bT(z)$ is in $(0,1)$ and $\bT(\rho)<\infty$.
He defines the bivariate function
\[
\bE(z,w) := z e^w\cdot \exp\Big(\sum_{m\ge 2}\bT(z^m)/m\Big),
\]
giving the recursion equation $\bT = \bE\big(z,\bT\big)$. Since
$\bE(z,w)$ is holomorphic in a neighborhood of $\bT$
he can invoke the Implicit Function Theorem to show that a necessary
condition for $z$ to be a dominant singularity, that is a singularity
on the circle of convergence, of $\bT$ is
\[
\bE_w\big(z,\bT(z)\big) = 1.
\]
From this P\'olya deduces that $\bT$ has a unique dominant singularity, namely
$z=\rho$. Next, since 
$\bE_z\big(\rho,\bT(\rho)\big),\, \bE_{ww}\big(\rho,\bT(\rho)\big)\neq 0$,
the Weierstra{\ss} Preparation Theorem shows that $\rho$ is a square-root
type singularity. Applying well known results derived from the 
Cauchy Integral Theorem 
\begin{equation}\label{cauchy}
t(n)\ =\ \frac{1}{2\pi i}\int_\cC \frac{\bT(z)}{z^{n+1}}dz
\end{equation}
one has the famous asymptotics 
\[
\begin{array}{l @{\quad} l}
\pmb{(\star)}&\qquad\qquad\qquad\qquad
t(n)\ \sim\ C \rho^{-n} n^{-3/2}\hspace*{2in}
\end{array}
\]
which occur so frequently in the study of recursively defined classes.

\subsection{Subsequent developments} 

Bender (\cite{Be}, 1974) proposed a general version of the P\'olya result, 
but Canfield (\cite{Canf}, 1983) found a flaw in the proof, and proposed a
more restricted version. Harary, Robinson and Schwenk (\cite{HRS}, 1975)  
gave a 20 step guideline on how to carry out a P\'olya style analysis of a 
recursion equation.
Meir and Moon (\cite{Me:Mo2}, 1989) made some further proposals on how to 
modify Bender's approach; in particular it was found that the hypothesis that 
{\em the coefficients of $\bE$ be nonnegative} was highly desirable, and covered 
a great number of important cases. 
This nonnegativity condition has continued to find
favor, being used in Odlyzko's survey paper \cite{Odly} and in the forthcoming book 
\cite{Fl:Se} of 
Flajolet and Sedgewick. Odlyzko's version seems to be a current standard---here it is 
(with minor corrections due to Flajolet and Sedgewick \cite{Fl:Se}). 

\begin{theorem}[Odlyzko \cite{Odly}, Theorem 10.6] 
\label{Odlyzko Thm}
Suppose 
\begin{eqnarray}
\bE(z,w)&=&\sum_{i,j\ge 0}e_{ij} z^i w^j \quad
\text{with } e_{00} = 0,\, e_{01} < 1,\, (\forall i,j)\, e_{ij}\ge 0\\ 
\bT(z)&=&\sum_{i\ge 1}t_i z^i\quad\text{with }(\forall i)\, t_i\ge 0
\end{eqnarray}
are such that
\begin{thlist}
\item $\bT(z)$ is analytic at $x=0$ 
\item $\bT(z) = \bE\big(z,\bT(z)\big)$
\item $\bE(z,w)$ is nonlinear in $w$
\item there are positive integers $i,j,k$ with $i<j<k$ such that
\begin{eqnarray*}
t_i,t_j, t_k&>&0\\
\gcd(j-i,k-i)&=&1.
\end{eqnarray*}
\end{thlist}
Suppose furthermore that there exist $\delta, r,s > 0$ such that
\begin{thlist}
\item [e] $\bE(z,w)$ is analytic in $|z|< r+\delta$ and $|w|<s+\delta$
\item [f] $\bE(r,s) = s$
\item [g] $\bE_w(r,s) = 1$
\item [h] $\bE_z(r,s) \neq 0$ and $\bE_{ww}(r,s)\neq 0$.
\end{thlist}
Then $r$ is the radius of convergence of $\bT$, $\bT(r) = s$, 
and as $n\rightarrow \infty$
\[
t_n\ \sim\ \sqrt{\frac{r \bE_z(r,s)}{2\pi 
\bE_{ww}(r,s)}} \cdot r^n n^{-3/2}.
\] 
\end{theorem}

\begin{remark}\label{Odly Rem}
As with P\'olya's original result, the asymptotics in these more general theorems
follow from information gathered on the location and nature of the dominant 
singularities of $\bT$. It has become popular to require that the solution $\bT$ 
have a {\bf unique} dominant singularity---to guarantee this happens the above 
theorem has the hypothesis \mbox{\rm(d)}. One can achieve this with a weaker hypothesis,
namely one only needs to require 
\[
\mbox{\rm(d$'$)} 
\quad \gcd\big(\{j-i : i< j\text{ and } t_i, t_j>0\}\big) = 1.
\]
Actually, given the other hypotheses of Theorem \ref{Odlyzko Thm}, the condition
\mbox{\rm(d$'$)} 
is necessary and sufficient that $\bT$ have a unique dominant singularity.
\end{remark}

The generalization of P\'olya's result that we find most convenient is given
in Theorem \ref{basic thm}.  We will also adopt the condition that $\bE$ 
have nonnegative coefficients, but point out that {\em under 
this hypothesis the 
location of the dominant singularities is quite easy to determine}. 
Consequently the unique singularity condition is not needed to determine the 
asymptotics.  

For further remarks on previous variations and 
generalizations of the work of P\'olya see $\S\,$\ref{history}.
The condition that the $\bE$ have nonnegative coefficients forces us to omit
the $\Set$ operator from our list of standard combinatorial operators. There
are a number of complications in trying to extend the results of this paper
to recursion equations $w = \bG(z,w)$ where $\bG$ has mixed signs appearing 
with its coefficients, including the problem of locating the
dominant singularities of the solution. The situation with mixed signs is
discussed in $\S$\,\ref{mixed}.

\subsection{Goal of this paper}

   Aside from the proof details that show we do not need to require that 
the solution $\bT$ have a unique dominant singularity, this paper is {\em not} 
about finding a
better way of generalizing P\'olya's theorem on trees. Rather the paper is concerned
with the ubiquity of the form $\pmb{(\star)}$ of asymptotics that 
P\'olya found for the recursively defined class of trees.\footnote{The motivation 
for our work came when a colleague, upon seeing the asymptotics of 
P\'olya for the first time, said ``Surely the form $\pmb{(\star)}$ hardly ever 
occurs! (when finding the asymptotics for the solution of an equation 
$w = \Theta(w)$ that recursively defines a class of trees)". A quick examination
of the literature, a few examples, and we were convinced that quite the
opposite held, that {\em almost any} reasonable class of trees defined by a 
recursive equation
that is nonlinear in $w$ would lead to an asymptotic law of P\'olya's form
$\pmb{(\star)}$.}

The goal of this paper is to exhibit a very large class of natural 
and easily recognizable operators $\Theta$ for 
which we can guarantee that a solution $w = \bT(z)$ to the recursion equation 
$w = \Theta(w)$ has coefficients that satisfy $\pmb{(\star)}$. By `easily
recognizable' we mean that you only have to look at the expression describing
$\Theta$---no further analysis is needed. This contrasts with the existing
literature where one is expected to carry out some calculations to determine
if the solution $\bT$ will have certain properties. For example, in Odlyzko's
version, Theorem \ref{Odlyzko Thm}, there is a great deal of work to be done, 
starting with checking that the solution $\bT$ is analytic at $z=0$.

In the formal specification theory for \emph{combinatorial classes} 
(see Flajolet and Sedgewick \cite{Fl:Se}) one starts with the binary operations 
of \emph{disjoint union} and \emph{disjoint sum} 
and adds unary
\emph{constructions} that transform a collection of objects (like trees) 
into a collection of objects (like forests). Such constructions
are \emph{admissible} if the generating series of the output class of the
construction is completely determined by the generating series of the 
input class.

We want to show that a recursive specification using almost any combination of these 
constructions, and others that we will introduce, yield classes whose
generating series have coefficients that obey the asymptotics
$\pmb{(\star)}$ 
of P\'olya. It is indeed a {\em universal law}.
The goal of this paper is to provide truly {\em practical} criteria 
(Theorem \ref{Main Thm})
to verify that many, if not most, of the common nonlinear recursion equations  lead to 
$\pmb{(\star)}$. 
Here is a contrived example to which this theorem applies:
\begin{equation}\label{ex1}
w\ =\  z \ +\ z \MSet\Big(\Seq\big(\sum_{n\in\Odd} 6^n w^n\big)\Big) \sum_{n\in\Even}(2^n+1)\big(\DCycle_{\Primes}(w)\big)^n\,.
\end{equation}
An easy application of Theorem \ref{Main Thm} 
(see $\S\,$\ref{applications}) 
tells us this particular recursion equation
has a recursively defined solution $\bT(z)$ with a positive radius 
of convergence, and the asymptotics for the coefficients $t_n$ have the form 
$\pmb{(\star)}$. 

The results of this paper apply to any combinatorial situation
described by a recursion equation of the type studied here. We put our 
focus on classes of trees because they are by far the most popular setting
for such equations.

\subsection{First definitions}

We start with our basic notation for number systems, power series and
open discs.

\begin{definition}
\begin{thlist}\itemsep=1ex
\item
$\bbR$ is the set of \underline{reals};
$\bbR^{\ge 0}$ is the set of \underline{nonnegative reals}. 
\item
$\bbP$ is the set of \underline{positive integers}. 
$\bbN$ is the set of \underline{nonnegative integers}. 
\item
$\bbR^{\ge 0}[[z]]$ is the set of 
power series in $z$ with nonnegative coefficients.
\item
$\rho_\bA$ is the \underline{radius {(of convergence)}}
of the power series $\bA$\,.
\item
For $\bA\in\bbR^{\ge 0}[[z]]$ we write 
$\bA\,=\,\sum_n a(n) z^n$ or $\bA\,=\,\sum_n a_n z^n$\,.
\item
For $r>0$ and $z_0\in\bbC$ the open disc of radius $r$ about $z_0$ is
$\bbD_r(z_0)\,:=\,\{z : |z-z_0|< r\}$ 
\end{thlist}
\end{definition}

\subsection{Selecting the domain }

We want to select a suitable collection of power series to work with when
determining solutions $w=\bT$ of recursion equations $w = \Phi(w)$.
The intended application is that $\bT$ be a generating series for some
collection of combinatorial objects.
Since generating series have nonnegative coefficients we naturally focus 
on series in $\bbR^{\ge 0}[[z]]$. 

There is one restriction that seems most desirable, namely to consider as generating
functions
only series whose constant term is 0. A generating series $\bT$ has
the coefficient $t(n)$ of $z^n$ counting (in some fashion) objects of size $n$.
It has become popular when working with combinatorial systems 
to admit a constant coefficient when it makes a result look simpler, for example
with permutations we write $\bA(z) = \exp\big(\bQ(z)\big)$, where $\bA(z)$ is the exponential
generating series for permutations, and $\bQ(z)$ the exponential generating series 
for cycles. 
$\bQ(z) = \log\big(1/(1-z)\big)$ will have a constant term 0, but $\bA(z)=1/(1-z)$ will have 
the constant term 1. 
Some authors like to introduce an `ideal' object of size 0 to go along with this
constant term.

There is a problem with this convention if one wants to look at compositions of
operators. For example, suppose you wanted to look at sequences of permutations.
The natural way to write the generating series would be to apply the sequence
operator $\Seq$ to $1/(1-z)$ above, giving $\sum 1/(1-z)^n$. 
Unfortunately this ``series'' has constant coefficient $=\infty$, so we 
do not have an analytical function. 
The culprit is the constant 1 in $\bA(z)$. If we drop the 1, so that we are 
counting only `genuine' permutations, the generating series for permutations 
is $z/(1-z)$; applying $\Seq$ to this gives $z/(1-2z)$, an analytical function 
with radius of convergence 1/2.

Consequently in this paper we return to the older convention of having the 
constant term be 0, so that we are only counting `genuine' objects.

\begin{definition}
For $\bA\in\bbR[[z]]$ we write $\bA\unrhd 0$ to say that all coefficients 
$a_i$ of $\bA$ are nonnegative. Likewise for $\bB\in\bbR[[z,w]]$ we write
$\bB\unrhd 0$ to say all coefficients $b_{ij}$ are nonnegative. Let
\begin{thlist}
\item
$\dom[z]\,:=\,\{\bA\in \bbR^{\ge 0}[[z]] : \bA(0) = 0\}$,
the set of power series $\bA \unrhd 0$ with constant term $0$; and let
\smallskip
\item
$\dom[z,w]\,:=\,\{\bE\in \bbR^{\ge 0}[[z,w]] : \bE(0,0) = 0\},$
the set of power series $\bE \unrhd 0$ with constant term $0$\,.
Members of this class are called \underline{elementary} power 
series.\footnote{We use the name {\em elementary} since a recursion 
equation
of the form $w = \bE(z,w)$ is in the proper form to employ the tools
of analysis that are presented in the next section.}
\end{thlist}
\end{definition}

When working with a member $\bE\in\dom[z,w]$
it will be convenient to use various series formats for writing $\bE$,
namely
\begin{eqnarray*}
\bE(z,w) &=& \sum_{ij} e_{ij} z^iw^j\\
\bE(z,w) &=& \sum_j \bE_j(z) w^j\\
\bE(z,w) &=& \sum_j \Big(\sum_i e_{ij} z^i \Big) w^j.
\end{eqnarray*}
This is permissible from a function-theoretic viewpoint 
since all coefficients $e_{ij}$ are 
nonnegative; for any given $z,w \ge 0$ the three formats converge to the
same value (possibly infinity).

An immediate advantage of working with series having nonnegative coefficients
is that the series is defined (possibly infinite) at its radius of convergence.

\begin{lemma}
For $\bT\,\in\,\dom[z]$ one has $\bT(\rho_\bT)\in[0,\infty]$.
Suppose $\bT(\rho_\bT)\,\in\,(0,\infty)$. Then
$\rho_\bT\,<\,\infty$\,; in particular $\bT$ is not a polynomial. 
If furthermore $\bT$ has integer coefficients then $\rho_\bT\,<\,1$.
\end{lemma}

\section{The theoretical foundations} \label{theor sect}

We want to show that the series $\bT$ that are recursively defined as 
solutions to functional equations $w = \bG(z,w)$ are such that
with remarkably frequency the
asymptotics of the coefficients $t_n$ are given by $\pmb{(\star)}$. 
Our main results deal with the
case that $\bG(z,w)$ is holomorphic in a neighborhood of $(0,0)$,
and the expansion $\sum g_{ij} z^iw^j$ is such that {\em all
coefficients $g_{ij}$ are nonnegative}. This covers most of the equations
arising from combinations of the popular combinatorial operators like 
Sequence, MultiSet and Cycle.

The referee noted that we had omitted one popular construction, namely
$\Set$, and the various restrictions $\Set_\bbM$ of $\Set$, and asked that 
we explain this omission. Although the equation $w = z+z \Set(w)$ has
been successfully analyzed in \cite{HRS}, there are difficulties when
one wishes to form composite operators involving $\Set$. These
difficulties arise from the fact that the resulting equation $w = \bG(z,w)$
has $\bG$ with coefficients having mixed signs. 
A general discussion of the mixed signs case is given in $\S$\,\ref{sub mixed}
and a particular discussion of the $\Set$ operator in $\S$\,\ref{sub set}.
Since the issue of 
mixed signs is so important we introduce the following abbreviations.

\begin{definition}
A bivariate series $\bE(z,w)$ and the associated functional equation
$w = \bE(z,w)$ are \underline{nonnegative} if the coefficients of 
$\bE$ are nonnegative. A bivariate series $\bG(z,w)$ and the associated
functional equation $w = \bG(z,w)$ have 
\underline{mixed signs} if some coefficients $g_{ij}$ are positive and 
some are negative.
\end{definition}

To be able to locate the difficulties when working with mixed signs, and
to set the stage for further research on this topic, we have put together
an essentially complete outline of the steps we use to prove that
a solution $\bT$ to a functional equation 
$w=\bE(z,w)$ satisfies the P\'olya asymptotics $\pmb{(\star)}$,
starting with the bedrock results of analysis such as the 
Weierstra{\ss} Preparation Theorem and the Cauchy Integral Formula. 
Although this background material has often been cited in work on 
recursive equations, it has never been written down in
a single unified comprehensive exposition. Our treatment of this background 
material goes beyond the existing literature to include a precise analysis of 
the nonnegative recursion equations whose solutions have multiple dominant 
singularities.

\subsection{A method to prove $\pmb{(\star)}$} 

Given $\bE\in\dom[z,w]$ and $\bT\in\dom[z]$ such that $\bT = \bE(z,\bT)$,
we use the following steps to show that the coefficients $t_n$
satisfy $\pmb{(\star)}$.
\begin{thlist}
\item
{\sc Show:} 
$\bT$ has radius of convergence $\rho := \rho_\bT > 0$.
\item
{\sc Show:} 
$\bT(\rho)<\infty$.
\item
{\sc Show:} 
$\rho <\infty$.
\item
{\sc Let:} $\bT(z) = z^d \bV(z^q)$ where $\bV(0) \neq 0$ and
$\gcd\big\{n : v(n) \neq 0\big\} = 1$.
\item
{\sc Let:} $\omega = \exp(2\pi i/q)$.
\item
{\sc Observe:} 
$\bT(\omega z) = \omega^d\bT(z)$, for $|z|<\rho$.
\item
{\sc Show:} 
The set of dominant singularities of $\bT$ is 
$\{z : z^q = \rho^q\}$.
\item
{\sc Show:} 
$\bT$ satisfies a quadratic equation, say 
\[
\bQ_0(z) + \bQ_1(z) \bT(z) + \bT(z)^2 \ =\ 0
\]
for $|z|<\rho$ and sufficiently near $\rho$, 
where $\bQ_0(z),\bQ_1(z)$ are analytic at $\rho$.
\item
{\sc Let:} $\bD(z) = \bQ_1(z)^2 - 4\bQ_0(z)$, the {\em discriminant} of 
the equation in (g).

\item
{\sc Show:} $\bD'(\rho)\neq 0$ in order to conclude that $\rho$ is a branch
point of order 2, that is, for $|z|<\rho$ and
sufficiently near $\rho$ one has
$\bT(z)\ =\ \bA(\rho - z) + \bB(\rho - z)\sqrt{\rho -z}$, 
where $\bA$ and $\bB$ are analytic at $0$, and 
$\bB(0)< 0$.
\item
{\sc Design:} A contour that is invariant under multiplication by $\omega$
to be used in the Cauchy Integral Formula to calculate $t(n)$.
\item
{\sc Show:} The full contour integral for $t(n)$ reduces to evaluating 
the portion lying between the angles $-\pi/q$ and $\pi/q$.
\item
{\sc Optional:}  One has a Darboux expansion for the asymptotics of $t(n)$.
\end{thlist}

Given that $\bE$ has nonnegative coefficients, items (a)--(f) can be
easily established by imposing modest conditions on $\bE$ (see 
Theorem \ref{basic thm}).
For (g) the method is to show that one has a functional
equation 
$\bF\big(z,\bT(z)\big)=0$
holding for $|z|\le\rho$ and sufficiently near $\rho$, that
$\bF(z,w)$ is holomorphic in a neighborhood of 
$\big(\rho,\bT(\rho)\big)$, and that
$\bF\big(\rho,\bT(\rho)\big)\,=\,
\bF_w\big(\rho,\bT(\rho)\big)\,=\,0$, but
$\bF_{ww}\big(\rho,\bT(\rho)\big)\,\neq\,0$.
These hypotheses allow one to apply the Weierstra{\ss} 
Preparation Theorem to obtain a quadratic equation for $\bT(z)$.

\begin{theorem}[Weierstra{\ss} Preparation Theorem]\label {WPThm}
Suppose $\bF(z,w)$ is 
a function of two complex variables and $(z_0,w_0)$
is a point in $\bbC^{\,2}$ such that:
\begin{thlist}
\item
$\bF(z,w)$ is holomorphic in a neighborhood of $(z_0,w_0)$
\item
$\displaystyle \bF(z_0,w_0)
\,=\, \frac{\partial \bF}{\partial w}(z_0,w_0)
\,=\,\cdots\,=\, \frac{\partial^{k-1} \bF}{\partial w^{k-1}}(z_0,w_0)\,=\,0$ 
\item
$\displaystyle 
\frac{\partial^k \bF}{\partial w^k}(z_0,w_0)\,\neq\,0$. 
\end{thlist}
Then in a neighborhood of $(z_0,w_0)$ one has $\bF(z,w)\,=\,\bP(z,w)\bR(z,w)$, 
a product of two holomorphic functions $\bP(z,w)$ and $\bR(z,w)$ where
\begin{itemize}
\item [(i)] $\bR(z,w)\neq 0$ in this neighborhood,
\item [(ii)] $\bP(z,w)$ is a `monic polynomial of degree $k$' in $w$, that is 
$\bP(z,w) \,=\, \bQ_0(z) + \bQ_1(z) w +\cdots + \bQ_{k-1}(z) w^{k-1} + w^k$,
and the $\bQ_i(z)$ are analytic in a neighborhood of $z_0$.
\end{itemize}
\end{theorem}
\begin{proof}
An excellent reference is Markushevich \cite{Marku}, Section 16, p.~105,
where one finds a leisurely and complete proof of the 
Weierstra{\ss} Preparation Theorem.
\end{proof}

There are two specializations of this result that we will be particularly
interested in: $k=1$ gives the Implicit Function Theorem, the best known
corollary of the Weierstra{\ss} Preparation Theorem;
and $k=2$ gives a quadratic equation for $\bT(z)$. 

\subsection{$k=1$: The implicit function theorem}

\begin{corollary} [k=1: Implicit Function Theorem] \label{IFThm}
Suppose $\bF(z,w)$ is 
a function of two complex variables and $(z_0,w_0)$
is a point in $\bbC^{\,2}$ such that:
\begin{thlist}
\item
$\bF(z,w)$ is holomorphic in a neighborhood of $(z_0,w_0)$
\item
$\displaystyle \bF(z_0,w_0)\,=\,0$
\item
$\displaystyle 
\frac{\partial \bF}{\partial w}(z_0,w_0)\,\neq\,0$.
\end{thlist}
Then there is an $\varepsilon>0$ 
and a function $\bA(z)$ such that for 
$z\,\in\,\bbD_\varepsilon(z_0)$, 
\begin{itemize}
\item[(i)]
$\bA(z)$ is analytic 
in $\bbD_\varepsilon(z_0)$ , 
\item[(ii)]
$\bF\big(z,\bA(z)\big)\,=\,0$
for $z\,\in\,\bbD_\varepsilon(z_0)$ , 
\item[(iii)]
for all $(z,w)\in \bbD_\varepsilon(z_0)\times \bbD_\varepsilon(w_0)$, 
if $\bF(z,w)\,=\,0$ then $w\,=\,\bA(z)$.
\end{itemize}
\end{corollary}
\begin{proof}
From Theorem \ref{WPThm} there is an $\varepsilon>0$ 
and a factorization of 
$\bF(z,w)\,=\,\bL(z,w)\bR(z,w)$, 
valid in 
$\bbD_\varepsilon(z_0)\times \bbD_\varepsilon(w_0)$, 
such that $\bR(z,w)\neq 0$ for
$(z,w)\,\in\,\bbD_\varepsilon(z_0)\times \bbD_\varepsilon(w_0)$, and
$\bL(z,w)\,=\,\bL_0(z)+ w$, with $\bL_0(z)$ analytic in 
$\bbD_\varepsilon(z_0)$. 

Thus $\bA(z)\,=\,-\bL_0(z)$ is such that 
$\bL\big(z,\bA(z)\big)\,=\,0$ on $\bbD_\varepsilon(z_0)$; so
$\bF\big(z,\bA(z)\big)\,=\,0$ on $\bbD_\varepsilon(z_0)$.
Furthermore, if $\bF(z,w)\,=\, 0$ with 
$(z,w)\in \bbD_\varepsilon(z_0)\times \bbD_\varepsilon(w_0)$, 
then $\bL(z,w)\,=\,0$, so $w\,=\,\bA(z)$.
\end{proof}

\subsection{$k=2$: The quadratic functional equation}

The fact that $\rho$ is an order 2 branch point comes
out of the $k=2$ case in the Weierstra{\ss} Preparation Theorem.

\begin{corollary}[$k=2$]\label{WPThm2}
Suppose $\bF(z,w)$ is 
a function of two complex variables and $(z_0,w_0)$
is a point in $\bbC^{\,2}$ such that:
\begin{thlist}
\item
$\bF(z,w)$ is holomorphic in a neighborhood of $(z_0,w_0)$
\item
$\displaystyle \bF(z_0,w_0) \,=\, \frac{\partial \bF}{\partial w}(z_0,w_0) \,=\,0$ 
\item
$\displaystyle 
\frac{\partial^2 \bF}{\partial w^2}(z_0,w_0)\,\neq\,0$. 
\end{thlist}
Then in a neighborhood of $(z_0,w_0)$ one has $\bF(z,w)\,=\,\bQ(z,w)\bR(z,w)$, 
a product of two holomorphic functions $\bQ(z,w)$ and $\bR(z,w)$ where 
\begin{itemize}
\item [(i)] $\bR(z,w)\neq 0$ in this neighborhood,
\item [(ii)] $\bQ(z,w)$ is a `monic quadratic polynomial' in $w$, that is 
$\bQ(z,w) \,=\, \bQ_0(z) + \bQ_1(z) w  + w^2$,
where $\bQ_0$ and $\bQ_1$ are analytic in a neighborhood of $z_0$.
\end{itemize}
\end{corollary}

\subsection{Analyzing the quadratic factor $\bQ(z,w)$} 

 \label{implicit}

Simple calculations are known 
(see \cite{Pl:Ro})
for finding all the partial derivatives of $\bQ$ and $\bR$ at
$\big(z_0,w_0\big)$
in terms of the partial derivatives of $\bF$ at the same point.
From this we can obtain important 
information about the coefficients of the discriminant $\bD(z)$
of $\bQ(z,w)$.

\begin{lemma}\label {WPThm2a}
Given the hypotheses (a)-(c) of Corollary \ref{WPThm2}
let 
$\bQ(z,w)$ and $\bR(z,w)$ be as described in (i)-(ii) of that corollary.
Then
\begin{itemize}
\item [(i)] $\bQ(z_0,w_0)\, =\, \bQ_w(z_0,w_0)\, =\, 0$
\item [(ii)] $\bR(z_0,w_0)\, =\, \bF_{ww}(z_0,w_0)/2$.
\end{itemize}
Let $\bD(z) \ =\ \bQ_1(z)^2 - 4\bQ_0(z)$, the discriminant of $\bQ(z,w)$.\\
Then
\begin{itemize}
\item [(iii)] $\bD(z_0) = 0$
\item [(iv)] $\bD'(z_0)\, =\, -8 \bF_z(z_0,w_0)\big/\bF_{ww}(z_0,w_0)$.
\end{itemize}
\end{lemma}

\begin{proof}
For (i) use Corollary \eqref{WPThm2} (b), the fact that $\bR(z_0,w_0) \neq 0$, and
\begin{eqnarray*}
\bF(z_0,w_0)&=&\bQ(z_0,w_0) \bR(z_0,w_0)\\
\bF_w(z_0,w_0)
&=& \bQ_w(z_0,w_0) \bR(z_0,w_0)\ +\ 
\bQ(z_0,w_0) \bR_w(z_0,w_0)\\
&=& \bQ_w(z_0,w_0) \bR(z_0,w_0).
\end{eqnarray*}
For (ii), since $\bQ$ and $\bQ_w$ vanish and $\bQ_{ww}$
evaluates to $2$ at $(z_0,w_0)$,
\begin{eqnarray*}
\bF_{ww}(z_0,w_0)
&=& 2\bR(z_0,w_0).
\end{eqnarray*}
For (iii) we have from (i)
\begin{eqnarray*}
0& =&\bQ_0(z_0) + \bQ_1(z_0) w_0 + {w_0}^2\\
0& =&\bQ_1(z_0)+ 2w_0
\end{eqnarray*}
and thus
\begin{eqnarray}
\bQ_1(z_0)& =&-2w_0 \label{A1}\\
\bQ_0(z_0)& =&{w_0}^2.\label{A0}
\end{eqnarray}
From \eqref{A1} and \eqref{A0} we have 
\begin{eqnarray*}
\bD(z_0)&=& \bQ_1(z_0)^2 - 4\bQ_0(z_0)\ =\ 4{w_0}^2 - 4 {w_0}^2\ =\ 0,
\end{eqnarray*}
which is claim (iii).

For claim (iv) start with
\begin{eqnarray*}
\bF_z(z_0,w_0)
&=&\bQ_z(z_0,w_0) \bR(z_0,w_0)\\
&=& \big(\bQ_0'(z_0) + w_0 \bQ_1'(z_0) \big)\bR(z_0,w_0).
\end{eqnarray*}
From the definition of $\bD(z)$ and \eqref{A1} 
\begin{eqnarray*}
\bD'(z_0)
&=& 2\bQ_1(z_0) \bQ_1'(z_0) - 4\bQ_0'(z_0)\\
&=& -4\big( \bQ_0'(z_0) + w_0 \bQ_1'(z_0) \big),
\end{eqnarray*}
so
\[
-4\bF_z(z_0,w_0)\ =\ \bD'(z_0) \bR(z_0,w_0).
\]
Now use (ii) to finish the derivation of (iv).

\end{proof}

\subsection{A square-root continuation of $\bT(z)$ when $z$ is near $\rho$}
Let us combine the above information into a proposition about 
a solution to a functional equation.

\begin{proposition} \label{crucial1}
Suppose $\bT\in\dom[z]$ is such that
\begin{thlist}
\item
$\rho:=\rho_\bT\in(0,\infty)$
\item
$\bT(\rho)<\infty$
\end{thlist}
and $\bF(z,w)$ is a function of two complex variables such that:
\begin{thlist}
\item[c]there is an $\varepsilon > 0$ such that $\bF\big(z,\bT(z)\big) = 0$
for $|z| < \rho$ and $|z-\rho| < \varepsilon$
\item[d]
$\bF(z,w)$ is holomorphic in a neighborhood of $\big(\rho,\bT(\rho)\big)$
\item[e]
$\displaystyle \bF\big(\rho,\bT(\rho)\big) \,=\, 
\frac{\partial \bF}{\partial w}\big(\rho,\bT(\rho)\big) \,=\,0$ 
\item[f]
$\displaystyle 
\frac{\partial \bF}{\partial z}\big(\rho,\bT(\rho)\big)\cdot
\frac{\partial^2 \bF}{\partial w^2}\big(\rho,\bT(\rho)\big)\,>\,0$. 
\end{thlist}
Then there are functions $\bA(z), \bB(z)$ analytic at $0$ 
such that
\[
\bT(z)\ =\ \bA(\rho - z) \,+\, \bB(\rho - z)\sqrt{\rho -z}
\]
for $|z|<\rho$ and near $\rho$
(see Figure \ref{puiseux}), and
\[
\bB(0)\ =\ -\sqrt{\frac{2\bF_z\big(\rho,\bT(\rho)\big)}
{\bF_{ww}\big(\rho,\bT(\rho)\big)}}\ <\ 0.
\]
\begin{figure}[h]
\centerline{\psfig{figure=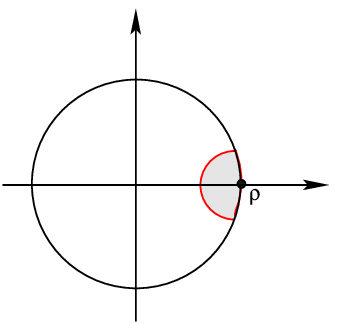}}
\caption{$\bT(z)\ =\ \bA(\rho-z) + \bB(\rho-z)\sqrt{\rho -z}$ in 
the shaded region}
\label{puiseux}
\end{figure}
\end{proposition}
\begin{proof}

Items (d)--(f) give the the hypotheses of Corollary \ref{WPThm2} 
with $(z_0,w_0) = \big(\rho,\bT(\rho)\big)$.
Let $\bQ_0(z)$, $\bQ_1(z)$ and $\bD(z)=\bQ_1(z)^2 - 4\bQ_0(z)$ 
be as in Corollary \ref{WPThm2}.
From conclusion (iv) of Lemma \ref{WPThm2a} we have 
\begin{equation}\label{D'}
\bD'(\rho)\ =\ -8\frac{\bF_z\big(\rho,\bT(\rho)\big)}
{\bF_{ww}\big(\rho,\bT(\rho)\big)}\ <\  0.
\end{equation}
From (c) and Corollary \ref{WPThm2}(i)
\[
\bQ_0(z)\,+\,\bQ_1(z)\bT(z)\,+\,\bT(z)^2\ =\ 0 
\]
holds in a neighborhood of $z=\rho$ intersected with 
${\bbD_\rho(0)}$
(as pictured in Figure \ref{puiseux}), so in this region
\[
\bT(z)\ =\ -\frac{1}{2}\bQ_1(z)\,+\,\frac{1}{2}\sqrt{\bD(z)}
\]
for a suitable branch of the square root.
Expanding $\bD(z)$ about $\rho$ gives
\begin{equation}\label{Dz}
\bD(z)\ =\ \sum_{k\ge 1} d_k(\rho -z)^k
\end{equation}
since $\bD(\rho)=0$ by (iii) of Lemma \ref{WPThm2a};
and $d_1 = -\bD'(\rho) > 0$ by \eqref{D'}.
Consequently 
\begin{equation}\label{sqrt cont}
\bT(z)\ =\ \underbrace{-\frac{1}{2}\bQ_1(z)}_{\displaystyle\bA(\rho - z)}\,
\underbrace{-\,\frac{1}{2}\sqrt{d_1}
\sqrt{1 + \sum_{k\ge 2}\frac{d_k}{d_1}(\rho - z)^{k-1}}}_{\displaystyle \bB(\rho-z)}
\,\cdot\,\sqrt{\rho -z}
\end{equation}
holds for $|z|<\rho$ and near $\rho$.  The negative sign of the second term is due to choosing
the branch of the square root which is consistent with the choice of branch implicit 
in Lemma \ref{binomial} when $\alpha = 1/2$, given that the $t(n)$'s are nonnegative.

Thus we have functions $\bA(z), \bB(z)$ analytic in a neighborhood
of $0$ with $\bB(0)\neq 0$ such that 
\[
\bT(z)\ =\ \bA(\rho-z) + \bB(\rho-z)\sqrt{\rho -z}
\]
for $|z|<\rho$ and near $\rho$. 
From \eqref{D'}, \eqref{Dz} and \eqref{sqrt cont}
\[
\bB(0)\ =\ -\frac{1}{2}\sqrt{d_1}\ =\ -\frac{1}{2}\sqrt{-\bD'(\rho)}\ =\  
-\sqrt{\frac{2\bF_z\big(\rho,\bT(\rho)\big)}
{\bF_{ww}\big(\rho,\bT(\rho)\big)}}\ <\  0.
\]
\end{proof}

Now we turn to recursion equations $w=\bE(z,w)$. 
So far in our discussion of the role of the Weierstra{\ss} Preparation Theorem we have not made any reference to the signs of the coefficients in the
recursion equation. The following proposition establishes a square-root 
singularity at $\rho$, and {\em the proof uses the fact that 
all coefficients of $\bE$ are nonnegative}. If we did not make
this assumption then items \eqref{Fww} and \eqref{Fz} below might fail to hold.
If \eqref{Fz} is false then $\bF_z\big(\rho,\bT(\rho)\big)$ may be
$0$, in which case
$\pmb{(\star)}$ fails. 
See section \ref{B0 condition} for a further discussion of this issue.

\begin{corollary} \label{crucial2}
Suppose $\bT\in\dom[z]$ and $\bE\in\dom[z,w]$
are such that
\begin{thlist}
\item
$\rho:=\rho_\bT\in(0,\infty)$
\item
$\bT(\rho)<\infty$
\item
$\bT(z) = \bE\big(z,\bT(z)\big)$ holds as an identity between formal
power series,
\item
$\bE(z,w)$ is not linear in $w$,
\item
$\bE_z\neq 0$
\item
$(\exists \varepsilon > 0)\,
\Big(\bE\big(\rho+\varepsilon,\bT(\rho)+\varepsilon\big)\,<\,\infty\Big)$.
\end{thlist}
Then there are functions $\bA(z), \bB(z)$ analytic at $0$ 
such that
\[
\bT(z)\ =\ \bA(\rho - z) + \bB(\rho - z)\sqrt{\rho -z}
\]
for $|z|<\rho$ and near $\rho$
(see Figure \ref{puiseux}), and
\[
\bB(0)\ =\ -\sqrt{\frac{2\bE_z\big(\rho,\bT(\rho)\big)}
{\bE_{ww}\big(\rho,\bT(\rho)\big)}}\ <\ 0.
\]
\end{corollary}
\begin{proof}

By (f) we can
choose $\varepsilon>0$ such that $\bE$ is holomorphic in
\[
\bbU\ =\  \bbD_{\rho + \varepsilon}(0) \times 
\bbD_{\bT(\rho) + \varepsilon}(0), 
\]
an open polydisc neighborhood of the graph of $\bT$.
Let 
\begin{equation}\label{def of F}
\bF(z,w)\ :=\ w \,-\,\bE(z,w).
\end{equation}
Then $\bF$ is holomorphic in $\bbU$, and one readily sees that
\begin{eqnarray}
\bF\big(z,\bT(z)\big)& =& \bT(z) \,-\,\bE\big(z,\bT(z)\big)\ =\ 0 
\quad\text{for }|z|\le \rho\label{F}\\
\bF_w(z,w)& =& 1 \,-\,\bE_w(z,w)\label{Fw}\\
\bF_{ww}\big(\rho,\bT(\rho)\big)& =& - 
\bE_{ww}\big(\rho,\bT(\rho)\big)\ <\ 0 
\quad\text{by (d) and $\bE\unrhd 0$} \label{Fww}\\
\bF_z\big(\rho,\bT(\rho)\big)&=& -\bE_z\big(\rho,\bT(\rho)\big)\ < \ 0 
\quad\text{by (e) and $\bE\unrhd 0$}\label{Fz}.
\end{eqnarray}
By Pringsheim's Theorem $\rho$ is a singularity of $\bT$. Thus
$\bF_w\big(\rho,\bT(\rho)\big)=0$ since one cannot use the Implicit Function 
Theorem to analytically continue
$\bT$ at $\rho$. 

We have satisfied the hypotheses of Proposition \ref{crucial1}---use 
\eqref{Fww} and \eqref{Fz} to obtain the formula for $\bB(0)$.

\end{proof}


\subsection{Linear recursion equations}

In a {\em linear} recursion equation
\[
w\ =\ \bA_0(z) + \bA_1(z) w
\]
one has
\begin{equation} \label{linear}
w\ =\ \frac{\bA_0(z)}{1 - \bA_1(z)}.
\end{equation}
From this we see that the collection of solutions to linear equations 
covers an enormous range.  For example, in the case
\[
w\ =\ \bA_0(z) + z w,
\]
any $\bT(z)\in\dom[z]$ with
nondecreasing eventually positive coefficients
is a solution to the above linear equation (which satisfies 
$\bA_0(z)+z w\,\unrhd\,0$) if we choose
$\bA_0(z) := (1-z) \bT(z)$.

When one moves to a $\Theta(w)$ that is
nonlinear in $w$, the range of solutions seems to be greatly constricted.
In particular with remarkable frequency one encounters solutions $\bT(z)$ 
whose coefficients are asymptotic to $C\rho^{-n}n^{-3/2}$.

\subsection{Binomial coefficients} \label{ultimate}

The asymptotics for the coefficients 
in the binomial expansion of $(\rho-z)^\alpha$ are the ultimate basis for 
the universal law $\pmb{(\star)}$. Of course if 
$\alpha\in\bbN$ then $(\rho-z)^\alpha$ is just a polynomial and the coefficients
are eventually 0. 

\begin{lemma}[See Wilf \cite{Wi}, p. 179] \label{binomial}
For $\alpha\in \bbR \setminus \bbN$ and $\rho\in(0,\infty)$
\[
[z^n]\,(\rho-z)^\alpha\ =\ 
(-1)^n
\binom{\alpha}{n}
\rho^{\alpha-n}\ 
\sim\ \frac{\rho^\alpha}{\Gamma(-\alpha)} \rho^{-n} n^{-\alpha-1}.
\]
\end{lemma}

\subsection{The Flajolet and Odlyzko singularity analysis}

In \cite{Fl:Od} Flajolet and Odlyzko develop {\em transfer} theorems via 
singularity analysis for functions $\bS(z)$ that 
have a unique dominant singularity.  
The goal is to develop a catalog of translations, or transfers,
 that say: {\em if $\bS(z)$ behaves like such and such 
near the singularity $\rho$
then the coefficients $s(n)$ have such and such asymptotic behaviour.}

\begin{figure}[ht]
\centerline{\psfig{figure=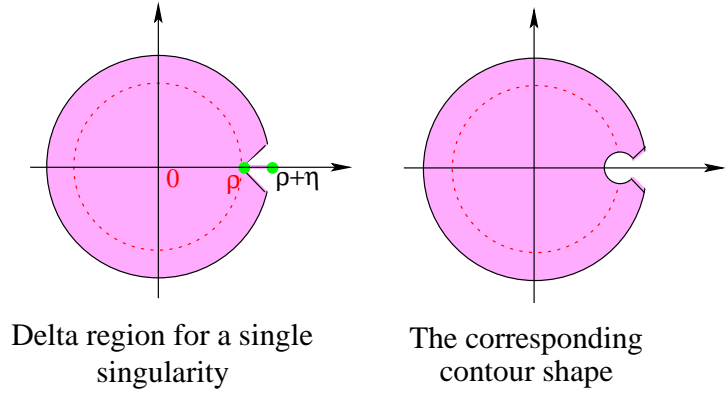}}
\caption{A Delta region and associated contour 
\label{delta}}
\end{figure}

Their work is based on applying the Cauchy Integral Formula to an
analytic continuation of $\bS(z)$ beyond its circle of convergence.
This leads to their basic notion of
a \underline{Delta neighborhood} $\Delta$ of $\rho$, that is,
a closed disc which is somewhat larger than the disc of radius $\rho$, 
but with an open pie shaped wedge cut out at the point $z=\rho$ 
(see Fig. \ref{delta}).  
We are particularly interested in their transfer theorem that directly
generalizes the binomial asymptotics given in Lemma \ref{binomial}. 

\begin{proposition}
[\cite{Fl:Od}, Corollary 2]
\label{FlOd prop}
Let $\rho\in(0,\infty)$ and suppose 
$\bS$ is analytic in $\Delta\setminus\{\rho\}$
where $\Delta$ is a Delta neighborhood of $\rho$.
If $\alpha\notin\bbN$ and
\begin{equation}
\label{alpha}
\bS(z)\ \sim \ K\big(\rho-z\big)^\alpha
\end{equation}
as $z\rightarrow \rho$ in $\Delta$, then
\[
s(n) \ \sim
[z^n]\,
K\big(\rho-z\big)^\alpha\ =\ 
(-1)^n K
\binom{\alpha}{n} 
\rho^{\alpha-n}\ 
\sim\ \frac{K\rho^\alpha}{\Gamma(-\alpha)}\cdot \rho^{-n}n^{-\alpha-1}.
\]
\end{proposition}

Let us apply this to the square-root singularities that
we are working with to see that one ends up with the 
asymptotics satisfying $\pmb{(\star)}$. 

\begin{corollary}
\label{FlOd cor}
Suppose $\bS\in\dom[z]$ has radius of convergence 
$\rho\in(0,\infty)$, and $\rho$ is the only dominant singularity 
of $\bS$. 
Furthermore suppose $\bA$ and $\bB$ are analytic at $0$ with 
$\bB(0) < 0$, $\bA(0)>0$ and 
\begin{equation}\label{sqrt}
\bS(z)\ = \bA(\rho -z) + \bB(\rho - z)\sqrt{\rho-z}
\end{equation}
for $z$ in some neighborhood of $\rho$, and $|z|<\rho$.

Then 
\[
s(n) \ \sim \ [z^n]\, \bB(0)\sqrt{\rho-z}
\ \sim\ \frac{-\bB(0) \sqrt{\rho}}{2\sqrt{\pi}} 
\cdot \rho^{-n} n^{-3/2}.
\]
\end{corollary}

\begin{proof}
One can find a Delta neighborhood $\Delta$ of $\rho$ 
(as in Fig. \ref{delta}) such that $\bS$ has an analytic 
continuation to $\Delta\setminus\{\rho\}$; and for 
$z\in\Delta$ and near $\rho$ one has \eqref{sqrt} holding.
Consequently
\[
\bS(z) - \bA(0)\ \sim\  \bB(0)\sqrt{\rho-z}
\]
as $z\rightarrow \rho$ in $\Delta$.
This means we can apply 
Proposition \ref{FlOd prop} to obtain
\begin{eqnarray*}
s(n) & \sim& \frac{\bB(0) \sqrt{\rho}}{\Gamma(-1/2)}
\cdot \rho^{-n} n^{-3/2}.
\end{eqnarray*}
\end{proof}

\subsection{On the condition $\bB(0)< 0$} \label{B0 condition}

In the previous corollary suppose that $\bB(0)= 0$ but $\bB\neq 0$.
Let $b_k$ be the first nonzero coefficient of $\bB$. 
The asymptotics for $s(n)$ are
\[
s(n)\ \sim\ b_k[z^n]\,\big(\rho - z\big)^{k +\frac{1}{2}},
\]
giving a law of the form $C \rho^{-n} n^{-k -\frac{3}{2}}$.
{\em We do not know of an example of $\bS$ defined by a nonlinear
functional equation that gives rise to such a solution with $k>0$, that is,
with the exponent of $n$ being $-5/2$, or $-7/2$, etc. Meir and Moon 
(p.~83 of \cite{Me:Mo2}, 1989) give the example 
\[
w\ =\ (1/6) e^w \sum_{n\ge 1} z^n/n^2
\]
where the solution $w=\bT$ has coefficient asymptotics given by
$t_n\sim C/n$.}

\subsection{Handling multiple dominant singularities}

We want to generalize Proposition \ref{FlOd prop} to cover the case of several dominant 
singularities equally spaced around the circle of convergence 
and with the function $\bS$ enjoying a certain kind of symmetry.

\begin{proposition} \label{FlOd Multi}
Given $q\in\bbP$ and $\rho\in(0,\infty)$ let 
\begin{eqnarray*}
\omega &:=&e^{2\pi i/q}\\
U_{q,\rho}& :=& \{ \omega^j\rho: j=0,1,\ldots,q-1\}.
\end{eqnarray*}
Suppose $\Delta$ is a generalized Delta-neighborhood of $\rho$
with wedges removed at the points in $U_{q,\rho}$ 
(see Fig. \ref{MultiDelta} for $q=3$),
\begin{figure}[h]
\centerline{\psfig{figure=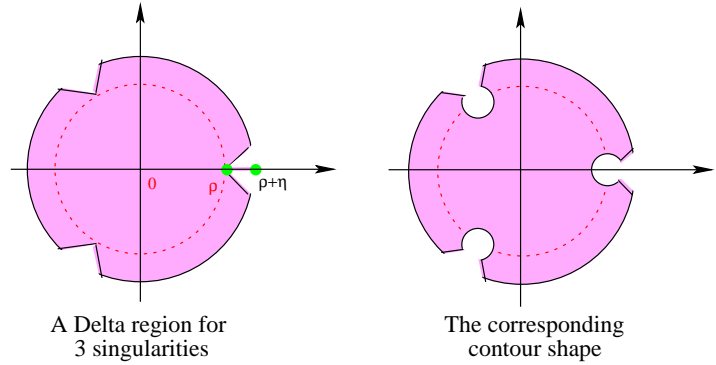}}
\caption{Multiple dominant singularities \label{MultiDelta}}
\end{figure}
 suppose
$\bS$ is continuous on $\Delta$ and
analytic in $\Delta\setminus U_{q,\rho}$, and suppose
$d$ is a nonnegative integer such that 
$
\bS\big(\omega z\big)\ = \omega^d \bS(z)
$
for $z\in\Delta$.

If
$\bS(z)\sim K(\rho-z)^\alpha$ as $z\rightarrow \rho$ in $\Delta$ 
and $\alpha\notin \bbN$
then 
\[
s(n) \ \sim\ \frac{q K \rho^\alpha}
{\Gamma(-\alpha)} \cdot \rho^{-n}
n^{-\alpha-1}\quad
\text{if } n\equiv d\mod q,
\]
$s(n)=0$ otherwise.
\end{proposition}

\begin{proof}
Given $\varepsilon > 0$ choose the contour $\cC$ to follow the 
boundary of $\Delta$ except for 
a radius $\varepsilon$ circular detour around each 
singularity $\omega^j\rho$ (see Fig. \ref{MultiDelta}). Then
\[
s(n)\ =\ \frac{1}{2\pi i}\int_{\cC} \frac{\bS(z)}{z^{n+1}}dz.
\]

Subdivide $\cC$ into $q$ congruent pieces 
$\cC_0,\ldots,\cC_{q-1}$ with $\cC_j$ centered around 
$\omega^j\rho$, choosing as 
the dividing points on $\cC$ the bisecting
points between successive singularities (see Fig. \ref{bisect2} for $q=3$).
\begin{figure}[ht]
\centerline{\psfig{figure=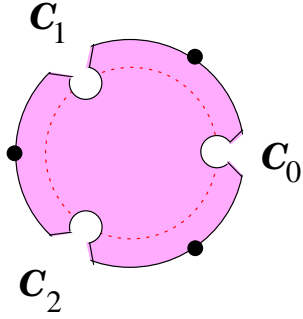}}
\caption{The congruent contour segments $\cC_j$\label{bisect2}}
\end{figure}
Then $\cC_j = \omega^j \cC_0$.
Let $s_j(n)$ be the portion of the integral for $s(n)$
taken over $\cC_j$, that is:
\[
s_j(n)\ =\ \frac{1}{2\pi i}\int_{\cC_j} \frac{\bS(z)}{z^{n+1}}dz.
\]
Then from $\bS(\omega z) = \omega^d \bS(z)$ and 
$\cC_j = \omega^j \cC_0$ we have
\begin{eqnarray*}
s_j(n)
& =& \frac{1}{2\pi i}\int_{\cC_j} \frac{\bS(z)}{z^{n+1}}dz\\
& =& \frac{1}{2\pi i}\int_{\cC_0} 
\frac{\omega^{dj} \bS(z)}{(\omega^j z)^{n+1}}\omega^j dz\\
& =& \omega^{j(d-n)} \frac{1}{2\pi i}\int_{\cC_0} 
\frac{\bS(z)}{z^{n+1}} dz\\
&=& \omega^{j(d-n)} s_0(n),
\end{eqnarray*}
so 
\begin{eqnarray*}
s(n)
& =&  \sum_{j=0}^{q-1} s_j(n)\\
& =&  
\Big(\sum_{j=0}^{q-1} \omega^{j(d-n)}\Big) s_0(n)\\
 & =&  \begin{cases} q s_0(n)&\text{if }n\equiv d \mod q\\
                     0&\text{otherwise.}
        \end{cases}
\end{eqnarray*}
We have reduced the integral calculation to the integral over
$\cC_0$, and this proceeds exactly as in \cite{Fl:Od} in the unique
singularity case described in Proposition \ref{FlOd prop}.
\end{proof}

Let us apply this result to the case of $\bS(z)$ having multiple
dominant singularities, equally spaced on the circle of convergence,
with a square-root singularity at $\rho$.

\begin{corollary} \label{MultiTransfer}

Given $q\in\bbP$ and $\rho\in(0,\infty)$ let 
\begin{eqnarray*}
\omega &:=&e^{2\pi i/q}\\
U_{q,\rho}& :=& \{ \omega^j\rho: j=0,1,\ldots,q-1\}.
\end{eqnarray*}
Suppose $\bS\in\dom[z]$ has radius of convergence 
$\rho\in(0,\infty)$, $U_{q,\rho}$ is the set of dominant 
singularities of $\bS$, and
$\bS(\omega z) = \omega^d \bS(z)$ for $|z|<\rho$ and 
for some $d\in\bbN$.

Furthermore suppose $\bA$ and $\bB$ are analytic at $0$ with 
$\bB(0) < 0$, $\bA(0) > 0$ and 
\begin{equation}\label{sqrt2}
\bS(z)\ = \bA(\rho -z) + \bB(\rho - z)\sqrt{\rho-z}
\end{equation}
for $z$ in some neighborhood of $\rho$, and $|z|<\rho$.
Then 
\begin{equation}\label{multi sqrt}
s(n) \  \sim  \ 
\frac{q \bB(0) \sqrt{\rho}}
{\Gamma(-1/2)} \cdot \rho^{-n}
n^{-3/2}\quad
\quad\text{for } n\equiv d \mod q.
\end{equation}
Otherwise $s(n)=0$.
\end{corollary}

\begin{proof}
Since the set of dominant singularities $U_{q,\rho}$ is finite
one can find a generalized Delta neighborhood $\Delta$ of $\rho$ 
(as in Fig. \ref{MultiDelta}) such that $\bS$ has a continuous extension
to $\Delta$ which is an analytic continuation to $\Delta\setminus U_{q,\rho}$; 
and for 
$z\in\Delta$ and near $\rho$ one has \eqref{sqrt2} holding.
Consequently
\[
\bS(z) - \bA(0) \ \sim\  \bB(0)\sqrt{\rho-z}
\]
as $z\rightarrow \rho$ in $\Delta$.
This means we can apply Proposition \ref{FlOd Multi} 
to obtain \eqref{multi sqrt}.
\end{proof}

\subsection{Darboux's expansion}

In 1878 Darboux \cite{darboux} published a procedure for expressing the
asymptotics of the coefficients $s(n)$ of a power series $\bS$ with
algebraic dominant singularities. 
Let us focus first on the case that $\bS$ has a single dominant
singularity, namely $z=\rho$, and it is of square-root type, say
\[
\bS(z)\ =\ \bA(\rho - z) + \bB(\rho - z) \sqrt{\rho - z}
\]
for $|z|< \rho$ and sufficiently close to $\rho$, where $\bA$ and $\bB$
are analytic at 0 and $\bB(0)< 0$.
From Proposition \ref{FlOd prop} we know that  
\[
s(n)\ =\ \big(1+\ro(1)\big) b(0) [z^n]\,\sqrt{\rho - z}.
\]

Rewriting the expression for $\bS(z)$ as
\[
\bS(z)\ =\ 
\sum_{j=0}^\infty \Big(a_j(\rho - z)^j \,+ \,
 b_j(\rho - z)^{j+\frac{1}{2}}\Big)
\]
we can see that the $m$th derivative of $\bS$ `blows up' as $z$ approaches
$\rho$ because the $m$th derivative of the terms on the right with $j< m$
involve terms with $\rho-z$ to a negative power. However for $j\ge m$ the
terms on the right have $m$th derivatives that behave nicely near $\rho$.
By shifting the troublesome terms to the left side of the equation, giving
\begin{eqnarray*}
\bS_m(z)& :=& \bS(z) -
\sum_{j<m} \Big(a_j(\rho - z)^j \,+ \,
 b_j(\rho - z)^{j+\frac{1}{2}}\Big)\\
& =& 
\sum_{j\ge m} \Big(a_j(\rho - z)^j \,+ \,
 b_j(\rho - z)^{j+\frac{1}{2}}\Big),
\end{eqnarray*}
one can see by looking at the right side that the $m$th 
derivative 
$\bS_m^{(m)}(z)$ 
of $\bS_m(z)$ has a square-root singularity at $\rho$ provided some
$b_j\neq 0$ for $j\ge m$.
Indeed 
$\bS_m^{(m)}(z)$ is very much like $\bS(z)$, being
analytic for $|z|\le \rho$ provided $z\neq \rho$.
If $b_m\neq 0$ we can apply Proposition \ref{FlOd prop} 
to obtain (for suitable $C_m$)
\[
[z^n]\, \bS_m^{(m)}(z) \ \sim\ C_m\rho^{-n} n^{-\frac{3}{2}}
\]
and thus
\[
[z^n]\, \bS_m(z) \ \sim\ C_m\rho^{-n} 
n^{-m-\frac{3}{2}}.
\]
This tells us that
\[
s(n)\ =\ \sum_{j<m} [z^n]\,\Big(a_j(\rho - z)^j \,+ \,
 b_j(\rho - z)^{j+\frac{1}{2}}\Big)\ +\ 
 \big(1+\ro(1)\big) C_m\rho^{-n} n^{-m-\frac{3}{2}}.
\]
For $n\ge m$ the part with the $a_j$ drops out, so we have
the Darboux expansion
\[
s(n)\ =\ \sum_{j<m} [z^n]\,\Big(
 b_j(\rho - z)^{j+\frac{1}{2}}\Big)\ +\ 
 \big(1+\ro(1)\big) C_m\rho^{-n} n^{-m-\frac{3}{2}}.
\]
The case of multiple dominant singularities is handled as previously.
Here is the result for the general exponent $\alpha$.

\begin{proposition}[Multi Singularity Darboux Expansion]\,
Given $q\in\bbP$
let 
\begin{eqnarray*}
\omega &:=&e^{2\pi i/q}\\
U_{q,\rho}& :=& \{ \omega^j\rho: j=0,1,\ldots,q-1\}.
\end{eqnarray*}
Suppose we have a generalized Delta-neighborhood $\Delta$
with wedges removed at the points in $U_{q,\rho}$ (see Fig. \ref{MultiDelta}) 
and
$\bS$ is analytic in $\Delta\setminus U_{q,\rho}$. Furthermore suppose
$d$ is a nonnegative integer such that 
$
\bS\big(\omega z\big)\ = \omega^d \bS(z)
$
for $|z|<\rho$.

If
\[
\bS(z)\ =\  \bA(\rho-z) + \bB(\rho-z)(\rho-z)^\alpha
\]
for $|z|<\rho$ and in a neighborhood of $\rho$,
and $\alpha\notin \bbN$,
then given $m\in\bbN$ with $b_m\neq 0$ there is a $C_m\neq 0$
such that
for $n\equiv d \mod q$
\[
s(n)\ =\ q\sum_{j<m} [z^n]\,\Big(
 b_j(\rho - z)^{j+\frac{1}{2}}\Big)\ +\ 
 \big(1+\ro(1)\big) C_m\rho^{-n} n^{-\alpha-(m+1)}.
\]
\end{proposition}

\subsection{An alternative approach: 
reduction to the aperiodic case} \label{altern}

In the literature one finds references to the option of using the 
aperiodic reduction $\bV$ of $\bT$, that is, using 
$\bT(z) = z^d\bV(z^q)$ where $\bV(0)\neq 0$ and $\gcd\{n : v(n) \neq 0\}=1$.
$\bV$ has a unique dominant singularity at 
$\rho_\bV\,=\,{\rho_\bT}^q$, so the hope 
would be that one could use a well known result like
Theorem \ref{Odlyzko Thm} to prove that $\pmb{(\star)}$ holds for $v(n)$.
Then $t(nq + d) = v(n)$ gives the asymptotics for the coefficients of $\bT$.

One can indeed make the transition from $\bT = \bE(z,\bT)$ to a functional
equation $\bV = \bH(z,\bV)$, but it is not clear if the property that 
$\bE$ is holomorphic at the endpoint of the graph of $\bT$ implies 
$\bH$ is holomorphic at the endpoint of the graph of $\bV$. 
Instead of the property
\[
(\exists \varepsilon > 0)\,\Big(\bE\big(\rho+\varepsilon,
\bT(\rho)+\varepsilon\big)\,<\,\infty\Big)
\]
of $\bE$ used previously, a stronger version seems to be needed, namely:
\[
(\forall y > 0)\Big[\bE\big(\rho,y\big)\,<\,\infty \Rightarrow\  
(\exists \varepsilon > 0)\,\Big(\bE\big(\rho+\varepsilon,
y+\varepsilon\big)\,<\,\infty\Big)\Big].
\]

We chose the singularity analysis because it sufficed to require
the weaker condition that $\bE$ be holomorphic at $\big(\rho,\bT(\rho)\big)$,
and because the expression for the constant term in the asymptotics was 
far simpler that what we obtained through the use of $\bV = \bH(z,\bV)$.
Furthermore, in any attempt to extend the analysis of the asymptotics to
other cases of recursion of equations one would like to have the ultimate
foundations of the Weierstra{\ss} Preparation Theorem and the Cauchy Integral
Theorem to fall back on.

\section{The Dominant Singularities of $\bT(z)$}

The recursion equations $w = \bE(z,w)$ we consider will be such that 
the solution $w=\bT$ has a radius of convergence $\rho$ in $(0,\infty)$ and finitely 
many dominant singularities, that is finitely many singularities on
the circle of convergence.
In such cases the primary technique to find the asymptotics
for the coefficients $t(n)$ is to apply Cauchy's Integral
Theorem \eqref{cauchy}.
Experience suggests that properly 
designed contours $\cC$ will concentrate the value of the integral \eqref{cauchy} 
on small 
portions of the contour near the dominant singularities of $\bT$---consequently great value 
is placed on locating the dominant singularities of $\bT$.

\begin{definition} 
For $\bT\in\dom[z]$ with radius $\rho\in (0,\infty)$ let
$\DSing(\bT)$ be the set of dominant singularities of $\bT$, that is, 
the set of singularities on the circle of convergence of $\bT$.
\end{definition}

\subsection{The spectrum of a power series}

\begin{definition} For $\bA\in\dom[z]$ let the \underline{spectrum} $\Spec(\bA)$
of $\bA$ be the set of $n$ such that the $n$th coefficient $a(n)$ is not
zero.\footnote{  
In the 1950s the logician Scholz defined 
the {\em spectrum} of a first-order sentence $\varphi$ to be the set of 
sizes of the finite models of $\varphi$. For example if $\varphi$ is an 
axiom for fields, then the spectrum would be the set of powers of primes.
There are many papers on this topic: a famous open problem due to Asser 
asks if the collection of spectra of first-order sentences is closed 
under complementation. This turns out to be equivalent to an open 
question in complexity theory. The recent paper \cite{Fi:Ma} of Fischer
and Makowsky has an
excellent bibliography of 62 items on the subject of spectra.

For our purposes, if $\bA(z)$ is a generating series for a class $\cA$ of 
combinatorial objects then the set of sizes of the objects 
in $\cA$ is precisely $\Spec(\bA)$.
}
It will be convenient to denote $\Spec(\bA)$ simply by $A$, so we have
\[
A \ =\ \Spec(\bA)\ =\ \{n : a(n) \neq 0\}.
\]
\end{definition}

In our analysis of the dominant singularities of $\bT$ it will be most
convenient
to have a simple calculus to work with the spectra of power series.

\subsection{An algebra of sets}

The spectrum of a power series from $\dom[z]$  is a subset of 
positive integers; the calculus
we use has certain operations on the subsets of the nonnegative integers.

\begin{definition}
For $I,J\subseteq \bbN $ and $j,m\in\bbN$ let
\begin{eqnarray*}
I+J& :=& \big\{i+j : i\in I, j\in J\big\}\\
I-j& :=& \big\{i-j : i\in I, \big\}\quad\text{where }j\le\min(I)\\
m \cdot J& :=& \{m\cdot j : j\in J\}\quad\text{for }m\ge 1\\
0\odot J &:=& \{0\}\\
m\odot J& :=& \underbrace{J+\cdots +J}_{m-{\rm times}} \quad \text{for $m \ge 1$}\\
I\odot J& :=& \bigcup_{i\in I} i\odot J\\
m| J&\Leftrightarrow & \big(\forall j\in J\big)\,(m|j).
\end{eqnarray*}
\end{definition}

\subsection{The periodicity constants}

Periodicity plays an important role in determining the dominant singularities. 
For example the generating series $\bT(z)$ of planar (0,2)-binary trees, 
that is, planar trees where each node has 0 or 2 successors, is defined by 
\[
\bT(z)\ =\ z + z\bT(z)^2.
\]
It is clear that all such trees have odd size, so one has
\[
\bT(z)\ =\ \sum_{j=0}^\infty t(2j+1) z^{2j+1} 
\ =\ z \sum_{j=0}^\infty t(2j+1) (z^2)^j.
\]
This says we can write $\bT(z)$ in the form 
\[
\bT(z)\ =\ z \bV(z^2).
\]
From such considerations one finds that $\bT(z)$ has exactly
two dominant singularities, $\rho$ and $-\rho$. (The general result is
given in Lemma \ref{dom sing}.)

\begin{lemma}\label{periodic form}
For $\bA\in \dom[z]$ let
\[
p\ :=\ \gcd A\qquad d\ :=\ \min A\qquad q\ :=\ \gcd (A-d).
\]
Then there are $\bU(z)$ and $\bV(z)$ in $\bbR^{\ge 0}[[z]]$ such that
\begin{thlist}
\item
$\bA(z) \ =\ \bU\big(z^p\big)$\quad with $\gcd(U) = 1$
\item
$\bA(z) = z^d\bV\big(z^q\big)$ with $\bV(0)\neq 0$ and $\gcd(V)=1$. 
\end{thlist}
\end{lemma}

\begin{proof} (Straightforward.)
\end{proof}

\begin{definition}
With the notation of Lemma \ref{periodic form},
$\bU(z^p)$ is the \underline{purely periodic form} of $\bA(z)$; and
$z^d\bV\big(z^q\big)$ 
is the \underline{shift periodic form} of $\bA(z)$.
\end{definition}

The next lemma is quite important---it says that the 
$q$ equally spaced points on the circle of convergence are all 
dominant singularities of $\bT$. 
Our main results depend heavily on the fact that the 
equations we consider are such that these are the {\em only} dominant 
singularities of $\bT$.

\begin{lemma} \label{min sing}
Let $\bT\in\dom[z]$ have radius of convergence $\rho\in(0,\infty)$
and the shift periodic form $z^d \bV(z^q)$. 
Then 
\[
\{z : z^q = \rho^q\} \ \subseteq\ \DSing(\bT).
\]
\end{lemma}

\begin{proof}
Suppose ${z_0}^q = \rho^q$ and suppose 
$\bS(z)$ is an analytic continuation 
of $\bT(z)$ into a neighborhood $\bbD_\varepsilon(z_0)$ of $z_0$. 
Let $\omega := z_0/\rho$. Then $\omega^q = 1$.
The function
$\bS_0(z) := \bS(\omega z)/\omega^d$ is an analytic function 
on $\bbD_\varepsilon(\rho)$. 
For $z\in \bbD_\varepsilon(\rho)\cap \bbD_\rho(0)$ we have 
\[
\omega z\in \bbD_\varepsilon(z_0)\cap \bbD_\rho(0),
\]
so
\[
\bS_0(z) \ = \ \bS(\omega z)/\omega^d\ =\ \bT(\omega z)/\omega^d\ =\ \bT(z).
\]
This means $\bS_0(z)$ is an analytic continuation of $\bT(z)$ at $z=\rho$, 
contradicting Pringsheim's Theorem that $\rho$ is a dominant singularity.
\end{proof}

\subsection{Determining the shift periodic parameters from $\bE$} \label{det shift}

\begin{lemma} \label{periodic constants}
Suppose $\bT(z) = \bE\big(z,\bT(z)\big)$ is a formal recursion that defines 
$\bT \in \dom[z]$, where $\bE\in\dom[z,w]$. 
Let the shift periodic form of $\bT(z)$ 
be $z^d  \bV(z^q)$. Then
\begin{eqnarray*}
d &=&\min(T)\ =\ \min(E_0)\\
q &=&\gcd(T-d)\ =\ \gcd\bigcup_{n\ge 0} \Big(E_n + (n-1)d\Big).
\end{eqnarray*}
\end{lemma}
\begin{proof}
Since $\bT$ is recursively defined by 
\[
\bT(z) \ = \ \sum_{n\ge 0} \bE_n(z) \bT(z)^n
\]
one has the first nonzero coefficient of $\bT$  being the first nonzero 
coefficient of $\bE_0$, and thus $d = \min(T) = \min(E_0)$. 
It is easy to see that we also have $q=\gcd(T-d)$.

Next apply the spectrum operator to the above functional equation to obtain
the set equation
\[
T\ =\ \bigcup_{n\ge 0} E_n + n\odot T, 
\]
and thus
\[
T-d\  =\  \bigcup_{n\ge 0} \Big(E_n + (n-1)d + n\odot (T-d)\Big). 
\]
Since $q = \gcd(T-d)$ it follows that 
$q | r := \gcd \Big(\bigcup_n E_n + (n-1)d\Big)$.

To show that $r | q$, and hence that $r = q$, note that
\[
w  =\  \bigcup_{n\ge 0} \Big(E_n + (n-1)d + n\odot w\Big)
\]
is a recursion equation whose unique solution is $w=T-d$.
Furthermore we can find the solution $w$ by iteratively applying
the set operator
\[
\Theta(w)\ :=\ \bigcup_{n\ge 0} \Big(E_n + (n-1)d + n\odot w\Big)
\]
to \O, that is, 
\[
T-d\ =\ \lim_{n\rightarrow \infty} \Theta^n(\text{\O}).
\]
Clearly $r\,|\,\text{\O}$, and a simple induction shows that for 
every $n$ we have $r\,\big|\,\Theta^n(\text{\O})$.
Thus $r\,|\,(T-d)$, so $r\,|\,q$, giving $r = q$. 
This finishes the proof that $q$ is the gcd of the
set $\bigcup_n\Big( E_n + (n-1)d\Big)$.
\end{proof}

\subsection{Determination of the dominant singularities}

The following lemma completely determines the dominant singularities
of $\bT$.
\begin{lemma}\label{dom sing}
Suppose
\begin{thlist}
\item
$\bT\in\dom[z]$ has radius of convergence $\rho\in (0,\infty)$ 
with $\bT(\rho) < \infty$, and
\item
$\bT(z) = \bE\big(z,\bT(z)\big)$, where $\bE\in\dom[z,w]$ 
is nonlinear in $w$ and holomorphic on (the graph of) $\bT$. 
\end{thlist}
Let the shift periodic form of $\bT(z)$ be $z^d  \bV(z^q)$.
Then 
\[
\DSing(\bT)\ =\ \{z : z^q = \rho ^ q\}.
\]
\end{lemma}

\begin{proof}
By the usual application of the implicit function theorem,
if $z$ is a dominant singularity of $\bT$  then 
\begin{equation}\label{orig sing cond}
\bE_w\big(z,\bT(z)\big) = 1.
\end{equation}
As
$\rho$ is a dominant singularity we can replace \eqref{orig sing cond} by
\begin{equation}\label{new sing cond}
\bE_w\big(z,\bT(z)\big) \ =\ \bE_w\big(\rho,\bT(\rho)\big).
\end{equation}

Let $\bU(z^p)$ 
be the purely periodic form of $\bE_w\big(z,\bT(z)\big)$. 
As the coefficients of $\bE_w$ are nonnegative it follows that 
\eqref{new sing cond} implies
\[
\DSing(\bT)\ \subseteq\ \{z : z^p = \rho^p\}.
\]
We know from Lemma \ref{min sing} that
\[
\{z : z^q = \rho^q\}\  \subseteq\ \DSing(\bT),
\]
consequently $q|p$. 

To show that $p \leq q$ first note that if $m\in \bbN$ then
\[
\gcd(m+T)\,\big|\,q.
\]
For if $r=\gcd(m+T)$ then for any $n\in T$ we have
$r|(m+n)$ and $r|(m+d)$. Consequently $r|(n-d)$, so $r|(T-d)$,
and thus $r|q$.

Since
\[
\bU(z^p)\ =\ \bE_w\big(z,\bT(z)\big)\ =\ \sum_{n\ge 1} \bE_n(z)  n 
\bT(z)^{n-1},
\]
applying the spectrum operator gives
\begin{eqnarray*}
\Spec\big(\bU(z^p)\big) 
& =& \bigcup_{n\ge 1} E_n +  (n-1)\odot T.
\end{eqnarray*}
Choose $n\ge 2$ such that $E_n\neq \text{\O}$ and choose $a\in E_n$. Then
\begin{eqnarray*}
\Spec\big(\bU(z^p)\big) 
& \supseteq& E_n + (n-1)\odot T\\
& \supseteq& \big(a+ (n-2)d\big) + T,
\end{eqnarray*}
so 
taking the $\gcd$ of both sides gives
\begin{eqnarray*}
p
&=&\gcd \Spec\big(\bU(z^p)\big)\\
&\le& \gcd\Big(\big(a+ (n-2)d\big) + T\Big)\,\Big|\,q.
\end{eqnarray*}
With $p =q$ it follows that we have proved the dominant
singularities are as claimed.
\end{proof}

\subsection{Solutions that converge at the radius of convergence}

The equations $w=\Theta(w)$ that we are pursuing
will have a solution $\bT$ that converges at the finite and positive 
radius of convergence $\rho_\bT$.

\begin{definition}
Let 
\[
\dom^\star[z]\ :=\ \{\bT \in \dom[z] : \rho_\bT\in (0,\infty),\ \bT(\rho_\bT)<\infty\}.
\]
\end{definition}

\subsection{A basic theorem}

The next theorem summarizes what we need from the preceding 
discussions 
to show that $\bT = \bE(z,\bT)$ leads to 
$\pmb{(\star)}$ holding for 
the coefficients $t_n$  of $\bT$.

\begin{theorem} \label{basic thm}
Suppose $\bT\in\dom[z]$ and $\bE\in\dom[z,w]$
are such that
\begin{thlist}
\item
$\bT(z) = \bE\big(z,\bT(z)\big)$ holds as an identity between formal
power series
\item
$\bT\in\dom^\star[z]$
\item
$\bE(z,w)$ is nonlinear in $w$
\item
$\bE_z\neq 0$ 
\item
$(\exists \varepsilon > 0)\,\Big(\bE\big(\rho+\varepsilon,
\bT(\rho)+\varepsilon\big)\,<\,\infty\Big)$.
\end{thlist}
Then 
\[
t(n)\ \sim\ q\sqrt
{\frac{\rho\bE_z\big(\rho,\bT(\rho)\big)} 
{2\pi \bE_{ww}\big(\rho,\bT(\rho)\big)}} 
\rho^{-n} n^{-3/2}
\quad\text{for } n\equiv d \mod q. 
\]
Otherwise $t(n)=0$.
Thus $\pmb{(\star)}$ holds on $\{n : t(n) > 0\}$.
\end{theorem}
\begin{proof}
By Corollary \ref{crucial2},
Corollary \ref{MultiTransfer} and Lemma \ref{dom sing}.
\end{proof}

\section{Recursion Equations using Operators}

Throughout the theoretical section, $\S\,$\ref{theor sect}, we only considered
recursive equations based on elementary operators $\bE(z,w)$.
Now we want to expand beyond these to include recursions that are based 
on popular combinatorial constructions
used with classes of unlabelled structures. As an umbrella concept
to create these various recursions we introduce the notion of 
{\em operators} $\Theta$. 

Actually if one is only interested in working with classes of labelled 
structures then it seems that the recursive equations based on
elementary power series are all that one needs. 
However, when working with classes of unlabelled structures, 
the natural way of writing down an equation corresponding to a recursive
specification is in terms of combinatorial operators like $\MSet$ 
and $\Seq$.  The resulting equation $w=\Theta(w)$, if properly designed,
 will have a unique solution $\bT(z)$ whose coefficients are recursively 
defined, and this solution will likely be needed to construct the 
translation of $w=\Theta(w)$ to an elementary recursion 
$w=\bE(z,w)$, a translation that is needed in order to apply the 
theoretical machinery of $\S\,$\ref{theor sect}.


\subsection{Operators} \label{comb op}

The mappings on generating series corresponding to combinatorial 
constructions are called operators. But we want to go beyond the 
obvious and include complex combinations of elementary and 
combinatorial operators. For this purpose we introduce 
a very general definition of an operator. 

\begin{definition}
An \underline{operator} is a mapping 
$\Theta: \dom[z]\,\rightarrow\, \dom[z]$\,. 
\end{definition}

Note that operators $\Theta$ act on $\dom[z]$\,, the set of formal 
power series with {\em nonnegative coefficients} and {\em constant term 0}. 
As mentioned before, the constraint that the constant terms of the power 
series be 0 makes for an elegant theory because compositions of operators 
are always defined.

A primary concern, as in the original work of P\'olya,
is to be able to handle combinatorial operators $\Theta$ that, 
when acting on $\bT(z)$, introduce terms like $\bT(z^2), \bT(z^3)$ etc.
For such operators it is natural to use power series $\bT(z)$ with
integer coefficients as one is usually working in the context of
ordinary generating functions. In such cases one has
$\rho\le 1$ for the radius of convergence of $\bT$,
provided $\bT$ is not a polynomial.

\begin{definition}
An \underline{integral operator} is a mapping 
$\Theta: \idom[z]\,\rightarrow\, \idom[z]$\,,
where 
$\idom[z]\ :=\ \big\{\bA\in\bbN[[z]] : \bA(0) = 0\big\}$,
the set of power series with nonnegative integer coefficients
and constant term zero.
\end{definition}

\begin{remark}
Many of the lemmas, etc, that follow have both a version for general 
operators and a version for integral operators. We will usually just 
state and prove the general version, leaving the completely parallel 
integral version as a routine exercise.
\end{remark}

\subsection{The arithmetical operations on operators} \label{arith sec}

The operations of addition, multiplication, positive scalar
multiplication and composition are defined on the set of operators 
in the natural manner:
\begin{definition}
\begin{eqnarray*}
(\Theta_1 \,+\, \Theta_2)(\bT)& :=& \Theta_1(\bT) \,+ \,\Theta_2(\bT)\\
(\Theta_1 \,\cdot\, \Theta_2)(\bT)& :=& \Theta_1(\bT) \,\cdot \,\Theta_2(\bT)\\
(c\cdot \Theta)(\bT)&:=&c\cdot \Theta(\bT)\\
(\Theta_1 \,\circ\, \Theta_2)(\bT)& :=& \Theta_1\big(\Theta_2(\bT)\big)\,,
\end{eqnarray*}
where the operations on the right side are the operations of formal power series. 
A set of operators is \underline{closed} 
if it is closed under the four arithmetical operations. 
\end{definition}

Note that when working with integral operators the scalars should be positive 
integers.
The operation of addition corresponds to the construction {\em disjoint union} 
and the operation of product to the construction {\em disjoint sum}, for both 
the unlabelled and the labelled case. 
Clearly the set of all [integral] operators is closed. 

\subsection{Elementary operators}

In a most natural way we can think of elementary power 
series $\bE(z,w)$ as operators.

\begin{definition}
Given $\bE(z,w)\in\dom[z,w]$ let the associated \underline{elementary
operator} be given by
\[
\bE : \bT\mapsto\bE(z,\bT)\quad\text{for }\bT\in\dom.
\]
\end{definition}

Two particular kinds of elementary operators are as follows. 
\begin{definition}
Let $\bA\in\dom[z]$.
\begin{thlist}
\item
The \underline{constant operator} $\Theta_\bA$ is given by
$\Theta_\bA:\bT\mapsto \bA$ for $\bT\in\dom[z]$, and
\item
the \underline{simple operator} $\bA(w)$ maps $\bT\in\dom[z]$ to
the power series that is the formal expansion of 
\[
\sum_{n\ge 1}a_n \Big(\sum_{j\ge 1} t_j z^j\Big)^n\,.
\]
\end{thlist}
\end{definition}

\subsection{Open elementary operators}

\begin{definition}
Given $a,b>0$,
an elementary operator $\bE(z,w)$ is \underline{open at $(a,b)$} if
\[
(\exists \varepsilon > 0)\Big(\bE(a+\varepsilon,b+\varepsilon)<\infty\Big).
\]
$\bE$ is \underline{open} if it is open at any
$a,b>0$ for which $\bE(a,b)<\infty$.
\end{definition}

Eventually we will be wanting an elementary operator to 
be open at $\big(\rho,\bT(\rho)\big)$ in order 
to invoke the Weierstra{\ss} Preparation Theorem. 

\begin{lemma} \label{cs open}
Suppose $\bA\in\dom[z]$ and $a,b>0$. \\
The constant operator $\Theta_\bA$ 
\begin{thlist}
\item
is open at $(a,b)$ iff
$ a<\rho_\bA$;
\item
it is open iff $\rho_\bA>0 \ \Rightarrow\ \bA(\rho_\bA)=\infty$.
\end{thlist}
\noindent
The simple operator $\Theta_\bA$ 
\begin{thlist}
\item[c]
is open at $(a,b)$ iff
$ b<\rho_\bA$; 
\item[d]
it is open iff $\rho_\bA>0 \ \Rightarrow\ \bA(\rho_\bA)=\infty$.
\end{thlist}
\end{lemma}

\begin{proof}
$\Theta_\bA$ is open at $(a,b)$ iff for some $\varepsilon > 0$
we have $\bA(a+\varepsilon)<\infty$. This is clearly 
equivalent to $a<\rho_\bA$. 

Thus $\rho_\bA>0$ and $\bA(\rho_\bA)<\infty$ 
imply $\Theta_\bA$ is not open at $(\rho_\bA,b)$ for any $b>0$,
hence it is not open. Conversely if $\Theta_\bA$ is not open then
$\rho_\bA>0$ and $\bA(a)<\infty$ 
for some $a,b>0$, but 
$\bA(a+\varepsilon)=\infty$ 
for any $\varepsilon>0$. This implies $a = \rho_\bA$. 

The proof for the simple operator $\bA(w)$ is similar.
\end{proof}

\subsection{Operational closure of the set of open $\bE$}

\begin{lemma} \label{closure of open E}
Let $a,b>0$. 
\begin{thlist}
\item
The set of elementary operators open at $(a,b)$
is closed under the arithmetical 
operations of 
scalar multiplication, addition and multiplication.
If\, $\bE_2$ is open at $(a,b)$ and $\bE_1$ is 
open at $\big(a,\bE_2(a,b)\big)$ then 
$\bE_1\big(z,\bE_2(z,w)\big)$ is open at $(a,b)$.
\item
The set of open elementary operators 
is  closed. 
\end{thlist}
\end{lemma}

\begin{proof}
Let $c>0$ and let $\bE,\bE_1,\bE_2$ be elementary operators 
open at $(a,b)$. 
Then 
\small
\begin{eqnarray*}
(\exists \varepsilon>0)\,\bE\big(a+\varepsilon,b+\varepsilon\big)<\infty
&\Rightarrow&(\exists \varepsilon>0)\,(c\bE)\big(a+\varepsilon,b+\varepsilon\big)<\infty\\
(\exists \varepsilon_1>0)\,\bE_1\big(a+\varepsilon_1,b+\varepsilon_1\big)<\infty
&\text{and}&
(\exists \varepsilon_2>0)\,\bE_2\big(a+\varepsilon_2,b+\varepsilon_2\big)<\infty
\\
&\Rightarrow&(\exists \varepsilon>0)
\bE_i\big(a+\varepsilon,b+\varepsilon\big)<\infty\ 
\text{for}\ i=1,2\\ 
&\Rightarrow&(\exists \varepsilon>0)\,
\big(\bE_1 + \bE_2\big)(a+\varepsilon,b+\varepsilon)<\infty\\
(\exists \varepsilon_1>0)\,\bE_1\big(a+\varepsilon_1,b+\varepsilon_1\big)<\infty
&\text{and}&
(\exists \varepsilon_2>0)\,\bE_2\big(a+\varepsilon_2,b+\varepsilon_2\big)<\infty
\\
&\Rightarrow&(\exists \varepsilon>0)
\bE_i\big(a+\varepsilon,b+\varepsilon\big)<\infty\ 
\text{for}\ i=1,2\\ 
&\Rightarrow&(\exists \varepsilon>0)\,
\big(\bE_1  \bE_2\big)(a+\varepsilon,b+\varepsilon)<\infty.
\end{eqnarray*}
\normalsize
Now suppose $\bE_2$ is open at $(a,b)$ and $\bE_1$ is 
open at $\big(a,\bE_1(a,b)\big)$. Then 
\small
\begin{eqnarray*}
(\exists \varepsilon_2>0)\,\bE_2\big(a+\varepsilon_2,b+\varepsilon_2\big)<\infty
&\text{and}&
(\exists \varepsilon_1>0)\,
\bE_1\big(a+\varepsilon_1,\bE_2(a,b)+\varepsilon_1\big)<\infty\\
&\Rightarrow&
(\exists \varepsilon>0)\,
\bE_1\big(a+\varepsilon,
\bE_2(a+\varepsilon,b+\varepsilon)+\varepsilon\big)<\infty.
\end{eqnarray*}
\normalsize
This completes the proof for (a). Part (b) is proved similarly.
\end{proof}

The {\em base} operators that we will use as a starting point
 are the elementary operators $\bE$ and all possible restrictions 
$\Theta_\bbM$ of the standard operators $\Theta$ of combinatorics 
discussed below. 
More complex operators called {\em composite} operators will be 
fabricated from these base operators by using the familiar 
{\em arithmetical operations} of addition, multiplication, scalar multiplication
and composition discussed in $\S\,$\ref{arith sec}.

\subsection{The standard operators on $\dom[z]$} \label{std op sec}

Following the lead of Flajolet and Sedgewick \cite{Fl:Se} we adopt as our
\emph{standard operators} 
$\MSet$ (multiset), $\Cycle$ (undirected cycle), 
$\DCycle$ (directed cycle) and $\Seq$ (sequence),
corresponding to the constructions by the same names.\footnote{Flajolet and Sedgewick
also include
$\Set$ as a standard operator, but we will not do so since, as mentioned in the second paragraph of $\S$\,\ref{theor sect}, for a given $\bT$, the 
series $\bG(z,w)$ associated with $\Set(\bT)$ may very well not be elementary.
For a discussion of mixed sign equations see $\S\,$\ref{mixed}.
}
These operators have well known analytic expressions, for example, 
\[
\begin{array}{l l}
\text{unlabelled multiset operator}&
1+\MSet(\bT)\ =\ \exp\Big(\sum_{j\ge 1} \bT(z^j)/j\Big)\\
\text{labelled multiset operator}&
\wMSet(\bT)\ =\ \sum_{j\ge 1}\bT(z)^j/j! \ =\ e^{\bT(z)} -1
\end{array}
\]

\subsection{Restrictions of standard operators} \label{restriction}

Let $\bbM\subseteq\bbP$\,. 
(We will always assume $\bbM$ is nonempty.)
The \emph{$\bbM$-restriction} of
a standard construction $\Delta$ applied to a class of trees means that 
one only takes those forests 
in $\Delta(\cT)$ where the number of trees is in $\bbM$\,. 
Thus $\MSet_{\{2,3\}}(\cT)$ takes all multisets of
two or three trees from $\cT$\,. 

The P\'olya {\em cycle index polynomials} 
$\bZ(\sH,z_1,\ldots,z_m)$
are very convenient for expressing such operators; such a polynomial 
is connected with a permutation group $\sH$ acting on an $m$-element 
set (see Harary and Palmer \cite{H:P}, p.~35). For $\sigma\in\sH$
let $\sigma_j$ be the number of $j$-cycles in a decomposition of
$\sigma$ into disjoint cycles. Then 
\[
\bZ(\sH,z_1,\ldots,z_m)\ :=\ \frac{1}{|H|}\sum_{\sigma\in H}
\prod_{j=1}^m {z_j}^{\sigma_j}.
\]
The only groups we consider are the following:
\begin{thlist}
\item
$\sS_m$ is the {\em symmetric group} on $m$ elements,
\item
$\sD_m$ the {\em dihedral group} of order $2m$, 
\item
$\sC_m$ the {\em cyclic group} of order $m$, and
\item
${\sf Id}_m$ the one-element {\em identity group} on $m$ elements.
\end{thlist}

The \emph{$\bbM$-restrictions} of the standard operators are
each of the form $\Delta_\bbM\,:=\,\sum_{m\in\bbM} \Delta_m$ where
$\Delta\in \{\MSet,\DCycle,\Cycle,\Seq\}$ and
$\Delta_m$ is given by:
{
\small
\[
\begin{array}{l l| l l}  
{\sf operator}&{\sf unlabelled\ case}&{\sf operator}&{\sf labelled\ case}\\
\MSet_m(\bT)
& \bZ\big(\sS_m,\bT(z),\ldots,\bT(z^m)\big) 
&\wMSet_m(\bT)
&(1/m!)\bT(z)^m\\
\Cycle_m(\bT)
&  \bZ\big(\sD_m,\bT(z),\ldots,\bT(z^m)\big)
&\wCycle_m(\bT)
& (1/2m)\bT(z)^m\\
\DCycle_m(\bT)
&  \bZ\big(\sC_m,\bT(z),\ldots,\bT(z^m)\big)
&\wDCycle_m(\bT)
& (1/m)\bT(z)^m\\
\Seq_m(\bT)
& \bZ\big({\sf Id}_m,\bT(z),\ldots,\bT(z^m)\big)
&\wSeq_m(\bT)
& \bT(z)^m\\
\end{array}
\]
\normalsize
}
Note that the labelled version of $\Delta_m$ is just
the first term of the cycle index polynomial for the unlabelled 
version,
and the sequence operators are the same in both cases.
We write simply $\MSet$ for $\MSet_\bbM$ if $\bbM$ is $\bbP$, etc.

In the labelled case the standard operators (with restrictions) are 
simple operators, whereas in the unlabelled case only $\Delta_{\{1\}}$ and the $\Seq_\bbM$ are simple.
The other standard operators in the unlabelled case are not elementary
because of the presence of terms $\bT(z^j)$ with $j>1$ when $\bbM\neq\{1\}$.


\subsection{Examples of recursion equations}

Table \ref{std examp} gives the recursion equations
for the generating series 
of several well-known classes of trees.
\begin{table}[h]
 \small
\[
\begin{array}{l @{\quad} l}
\text{Recursion Equation}&\text{Class of Rooted Trees}\\ 
\hline
w\,=\,z\,+\,z w&\text{chains}\\
w\,=\,z\,+\,z \Seq(w)&\text{planar}\\
w\,=\,m z\,+\,m z \Seq(w)&\text{$m$-flagged 
planar\footnotemark}\\
w\,=\,z e^w&\text{labelled}\\
w\,=\,z\,+\,z \MSet(w)&\text{unlabelled }\\
w\,=\,z\,+\,z \MSet_{\{2,3\}}(w)&\text{unlabelled (0,2,3)-}\\
w\,=\,z\,+\,z\Seq_2(w)&\text{unlabelled binary planar}\\
w\,=\,z\,+\,z\MSet_2(w)&\text{unlabelled binary }\\
w\,=\,z\,+\,z w^2&\text{labelled binary }\\
w\,=\,z\,+\,z\big(w\,+\,\MSet_2(w)\big)&\text{unlabelled unary-binary }\\
w\,=\,z\,+\,z\MSet_r(w)&\text{unlabelled $r$-regular }\\
\hline
\end{array}
\]
\caption{Familiar examples of recursion equations\label{std examp}}
\end{table} \footnotetext{\emph{$m$-flagged} 
means one can attach any subset of $m$ given flags to each vertex.
This is just a colorful way of saying that the tree structures are augmented
with $m$-unary predicates $U_1,\ldots,U_m$\,, and each can hold on any subset
of a tree independently of where the others hold. }


\subsection{Key properties of operators}

Now we give a listing of the various properties of abstract 
operators that are needed to prove a universal law for recursion equations. 
The first question to be addressed is 
``Which properties does $\Theta$ need in order to guarantee that 
$w=\Theta(w)$ has a solution?''

\subsection{Retro operators} 

There is a simple natural property of an operator $\Theta$ that guarantees
an equation $w = \Theta(w)$ has a {\em unique} solution that is
determined by a recursive computation of the coefficients, namely
$\Theta$ calculates, given $\bT$, the $n$th coefficient of 
$\Theta(\bT)$ solely on the basis of the values of $t(1),\ldots,t(n-1)$. 

\begin{definition}
An operator $\Theta$ is 
\underline{retro}
if there is a sequence $\sigma$ of functions such that for 
$\bB\,=\,\Theta(\bA)$ one has 
$b_n\,=\,\sigma_n(a_1,\ldots,a_{n-1})$\,, where $\sigma_1$ is a
constant\,.
\end{definition}

There is a strong temptation to call such $\Theta$ \emph{recursion} 
operators since they will be used to recursively define
generating series. But without the context of a recursion equation
there is nothing recursive
about $b_n$ being a function of $a_1,\ldots,a_{n-1}$\,.

\begin{lemma}\label{retro lem}
A retro operator $\Theta$ has a unique fixpoint in $\dom[z]$, 
that is,
there is a unique power series $\bT\,\in\,\dom[z]$ such that
$
\bT\,=\, \Theta(\bT)\,.
$
We can obtain $\bT$ by an iterative application of $\Theta$
to the constant power series $0$:
\[
\bT\ = \ \lim_{n\rightarrow \infty} \Theta^n(0)\,.
\]
If $\Theta$ is an integral retro operator then $\bT\in\idom[z]$.
\end{lemma}

\begin{proof}
Let $\sigma$ be the sequence of functions that witness the
fact that $\Theta$ is retro.
If $\bT=\Theta(\bT)$ then 
\begin{eqnarray*}
t(1)&=&\sigma_1\\
t(n)&=&\sigma_n\big(t(1),\ldots,t(n-1)\big)\quad\text{for }n>1.
\end{eqnarray*}
Thus there is at most one possible fixpoint $\bT$ of $\Theta$; and 
these two equations show how to recursively find such a $\bT$.

A simple argument shows that $\Theta^{n+k}(0)$ agrees with $\Theta^n(0)$
on the first $n$ coefficients, for all $k\ge 0$\,. Thus 
$\lim_{n\rightarrow \infty} \Theta^n(0)$ is a fixpoint, 
and hence {\em the} fixpoint\,. If $\Theta$ is also integral then
each stage $\Theta^n(0)\in\idom[z]$, so $\bT\in\idom[z]$.
\end{proof}

Thus if $\Theta$ is a retro operator then the functional equation
$ w = \Theta(w)$ has a unique solution $\bT(z)$. 
Although the end goal is to have an equation $w=\Theta(w)$ with $\Theta$ a retro
operator, for the intermediate stages it is often more desirable to work with
{\em weakly} retro operators. 

\begin{definition}
An operator $\Theta$ is 
\underline{weakly retro}
if there is a sequence $\sigma$ of functions such that for 
$\bB\,=\,\Theta(\bA)$ one has 
$b_n\,=\,\sigma_n(a_1,\ldots,a_n)$\,.
\end{definition}

\begin{lemma} \label{retro is full}
\begin{thlist}
\item
The set of retro operators is closed.
\item
The set of weakly retro operators is  closed and
includes all elementary operators and all restrictions of standard operators.
\item
If $\Theta$ is a weakly retro operator then $z \Theta$ and 
$w\Theta$ are both retro operators.
\end{thlist}
\end{lemma}

\begin{proof}
For (a),
given retro operators $\Theta,\Theta_1,\Theta_2$, a positive constant $c$ and a
power series $\bT\unrhd 0$, we have
\begin{eqnarray*}
{[z^n]\,}(c \Theta)(\bT) &=&c\big([z^n]\,\Theta(\bT)\big)\\
{[z^n]\,}(\Theta_1+\Theta_2)(\bT)& =& [z^n]\,\Theta_1(\bT) + [z^n]\,\Theta_2(\bT)\\
{[z^n]\,}(\Theta_1\Theta_2)(\bT)
& =& \sum_{j=1}^{n-1}[z^j]\,\Theta_1(\bT)  [z^{n-j}]\Theta_2(\bT)\\
{[z^n]\,}(\Theta_1\circ\Theta_2)(\bT)
& =& \sigma_n
\big([z^1]\Theta_2(\bT),\ldots,[z^n]\,\Theta_2(\bT)\big),
\end{eqnarray*}
where $\sigma$ is the sequence of functions that witness the fact that
$\Theta_1$ is a retro operator.
In each case it is clear that the value of the right side depends only 
on the first $n-1$ coefficients of $\bT$. Thus the set of retro 
operators is  closed.
 
For (b) use the same proof as in (a), after changing the initial operators to
weakly retro operators, to show that the set of weakly retro operators is closed.

For an elementary operator $\bE(z,w)$ and power series $\bT\unrhd 0$
we have, after writing $\bE(z,w)$ as 
$\sum_{i\ge 0} \bE_i(z) w^i$,
\begin{eqnarray*}
{[z^n]\,}\bE\big(z,\bT(z)\big)
&=&[z^n]\, \sum_{j\ge 0} \bE_j(z) \bT(z)^j\\
&=&\sum_{j\ge 0} \sum_{i= 0}^n e_{ij} \big[z^{n-i}\big]\bT(z)^j.
\end{eqnarray*}
The last expression clearly depends only on the first $n$ 
coefficients of $\bT(z)$. 
Thus all elementary operators are weakly retro operators.

Let $Z(\bH,z_1,\ldots,z_m)$ be a cycle index polynomial. 
Then for $\bT\in\dom[z]$ one has 
\[
{[z^n}]\bT(z^j)\ =\ 
\begin{cases} 
0&\text{if $j$ does not divide $n$}\\
t(n/j)&\text{if }j|n.
\end{cases}
\]
Thus the operator that maps $\bT(z)$ to $\bT(z^j)$ is 
a weakly retro operator.
The set of weakly retro operators is closed,
so the operator mapping $\bT$ to 
$Z\big(\bH,\bT(z),\ldots,\bT(z^m)\big)$ is weakly retro.
Now every restriction $\Delta_\bbM$ of a standard operator is a
(possibly infinite) sum of such instances of cycle index polynomials;
thus they are also weakly retro.

For (c) note that 
\begin{eqnarray*}
{[z^n]}\,\big(z \Theta(\bT)\big)
&=&[z^{n-1}]\,\Theta(\bT)\\
{[z^n]}\,\big(\bT \Theta(\bT)\big)
&=&\sum_{j=1}^{n-1} t_j [z^{n-j}]\,\Theta(\bT),
\end{eqnarray*}
and in both cases the right side depends only on $t_1,\ldots,t_{n-1}$.

\end{proof}

\begin{lemma}\label{retro form}
\begin{thlist}
\item
An elementary operator $\bE(z,w)=\sum_{ij} e_{ij}z^iw^j$ is 
retro iff $e_{01} = 0$. 
\item
A restriction $\Delta_\bbM$ of a standard operator $\Delta$ is 
retro iff $1\notin\bbM$.
\end{thlist}
\end{lemma}
\begin{proof}
For (a) let $\bT\in\dom[z]$. Then
\[
[z^n]\,\bE\big(z,\bT(z)\big)\ =\ 
\sum_{j\ge 0} \sum_{i= 0}^n e_{ij} \big[z^{n-i}\big]\bT(z)^j,
\]
which does not depend on $t(n)$ iff $e_{01}=0$.

For (b) one only has to look at the definition of the P\'olya cycle index
polynomials. 
\end{proof}

The property of being retro for an elementary $\bE(z, w)$ is very closely related to the necessary and sufficient conditions for an equation $w = \bE(z,w)$ to give a recursive definition of a function $\bT\in\dom[z]$ that is not 0.  To see this rewrite the equation in the form
\begin{eqnarray*}
(1-e_{01}) w
& =& \big(e_{10}z + e_{20}z^2 + \cdots \big) +
\big(e_{11}z + \cdots \big) w + \cdots.
\end{eqnarray*}
(We know that $e_{00}=0$ as $\bE$ is elementary.)
So the first restriction needed on $\bE$ is that $e_{01}< 1$.

Suppose this condition on $e_{01}$ holds.
Dividing through by $1-e_{01}$ gives an {\em equivalent} equation
with no occurrence of the linear term $z^0w^1$ on the right hand side,
thus leading to the use of
$
e_{01} = 0
$
rather than the apparently weaker condition $e_{01} < 1$.

To guarantee a nonzero solution we also need that
$
\bE_0(z) \neq 0,
$
and by the recursive construction these conditions suffice.

Now that we have a condition, being retro, to guarantee that 
$w=\Theta(w)$ is a recursion equation with a unique solution 
$w=\bT$, the next goal is to find simple conditions on 
$\Theta$ that ensure this solution will have the desired 
asymptotics.


\subsection{Dominance between power series}
It is useful to have a notation to indicate that 
the coefficients of one series dominate those of another.

\begin{definition} 
For power series $\bA,\bB\,\in\,\dom[z]$ 
we say $\bB$ \underline{dominates} $\bA$\,,
written $\bA\,\unlhd\,\bB$\,, if $a_j\,\le\, b_j$ for all $j$.

Likewise for power series $\bG,\bH\,\in\,\dom[z,w]$ 
we say $\bH$ \underline{dominates} $\bG$\,,
written $\bG\,\unlhd\,\bH$\,, if $g_{ij}\,\le\, h_{ij}$ for all $i,j$.
\end{definition}

\begin{lemma} \label{dom series}
The dominance relation $\unlhd$ is a partial ordering on $\dom[z]$
preserved by the arithmetical operations: 
for $\bT_1,\bT_2,\bT\in\dom[z]$ and a constant $c>0$,
if 
$\bT_1\unlhd\bT_2$ then
\begin{eqnarray*}
c\cdot\bT_1&\unlhd& c\cdot\bT_2\\
\bT_1 + \bT& \unlhd& \bT_2 + \bT\\
\bT_1 \cdot \bT& \unlhd& \bT_2 \cdot \bT\\
\bT_1 \circ\bT& \unlhd& \bT_2 \circ \bT\\
\bT \circ\bT_1& \unlhd& \bT \circ \bT_2.
\end{eqnarray*}
\end{lemma}
\begin{proof}
Straightforward.
\end{proof}

\subsection{The dominance relation on the set of operators}

\begin{definition}  
\begin{thlist}
\item
For operators $\Theta_1,\Theta_2$ 
we say
$\Theta_2$ \underline{dominates} $\Theta_1$\,, 
symbolically $\Theta_1\,\ole\,\Theta_2$\,, 
if for any $\bT\in \dom[z]$ one has
$\Theta_1(\bT)\,\unlhd\,\Theta_2(\bT)$\,.
\item
For integral operators $\Theta_1,\Theta_2$ 
we say
$\Theta_2$ \underline{dominates} $\Theta_1$\,, 
symbolically $\Theta_1\,\iole\,\Theta_2$\,, 
if for any $\bT\in \idom[z]$ one has
$\Theta_1(\bT)\,\unlhd\,\Theta_2(\bT)$\,.
\end{thlist}
\end{definition}

As usual we continue our discussion mentioning only the general 
operators when the integral case is exactly parallel.
It is straightforward to check that the dominance relation $\ole$ is a 
partial ordering on the set of operators which is preserved by 
addition, multiplication and positive scalar multiplication.
Composition on the right also preserves $\ole$, that is, for operators
$\Theta_1 \ole \Theta_2$ and $\Theta$,
\[
\Theta_1 \circ\Theta \ \ole \ \Theta_2 \circ \Theta.
\]
However composition on the left requires an additional property, \emph{monotonicity}.

The bivariate $\bE$ in $\dom[z,w]$ play a dual role, on the one
hand simply as power series, and on the other as operators.
Each has a notion of dominance, and they are related.

\begin{lemma} \label{bidom}
For $\bE,\bF\in\dom[z,w]$ we have
\[
\bE \unlhd \bF\ \Rightarrow\ \bE \ole\bF.
\]
\end{lemma}
\begin{proof}
Suppose $\bE \unlhd \bF$ and let $\bT\in\dom[z]$.
Then
\[
[z^n]\,\bE(z,\bT)
\ =\ \sum e_{ij}[z^{n-i}]\bT(z)^j
\ \le\ \sum f_{ij}[z^{n-i}]\bT(z)^j
\ =\ [z^n]\,\bF(z,\bT),
\]
so $\bE(z,\bT)\unlhd\bF(z,\bT)$. As $\bT$ was arbitrary, $\bE\ole\bF$.
\end{proof}

\subsection{Monotone operators}

\begin{definition}
An operator $\Theta$ 
is \underline{monotone} 
if it preserves $\unlhd$\,, that is, $\bA\,\unlhd\,\bB$ implies 
$\Theta(\bA)\,\unlhd\,\Theta(\bB)$ for $\bA,\bB\in\dom[z]$\,. 
\end{definition}

\begin{lemma}
If $\Theta_1 \ole\Theta_2$ and $\Theta$ is monotone then
\[
\Theta \circ\Theta_1 \ \ole\ \Theta \circ \Theta_2.
\]
\end{lemma}

\begin{proof}
Straightforward.
\end{proof}

\begin{lemma} \label{monotone closed}
The set of monotone operators is  closed and includes
all elementary operators and all restrictions of the standard operators.
\end{lemma}

\begin{proof}
Straightforward.
\end{proof}

\subsection{Bounded series}

\begin{definition}\label{AR}
For $R>0$ let $\bA_R(z) := \sum_{n\ge 1}R^n\,z^n$.
A series $\bT\in\dom[z]$ is 
\underline{bounded} if $\bT\unlhd\bA_R$ for some $R>0$.
\end{definition}

An easy application of the Cauchy-Hadamard Theorem shows that $\bT$ is
bounded iff it is analytic at $0$.

The following basic facts about the series $\bA_R(z)$ show
that the collection of bounded series is closed under the 
arithmetical operations, a well known fact. Of more interest
will be the application of this to the collection of bounded
operators in Section \ref{bdd closed sect}.
\begin{lemma}\label{bounded updir}
For $c,R,R_1,R_2>0$
\begin{eqnarray*}
R_1\le R_2&\Rightarrow&\bA_{R_1}\unlhd\bA_{R_2}\\
c \bA_R& \unlhd& \bA_{(c+1)R}\\ 
\bA_{R_1} + \bA_{R_2}&\unlhd& \bA_{R_1+R_2}\\
\bA_{R_1}  \bA_{R_2}&\unlhd& \bA_{R_1+R_2 }\\
\bA_{R_1} \circ\bA_{R_2}&\unlhd& \bA_{ 2(1+R_1+R_2)^2}.
\end{eqnarray*}
\end{lemma}
\begin{proof}
The details are quite straightforward---we give the proofs for the
last two items.
\begin{eqnarray*}
\big(\bA_{R_1}  \bA_{R_2}\big)(z) 
&=&
 \Big(\sum_{j\ge 1}{R_1}^j z^j\Big)\ 
\cdot\  \Big(\sum_{j\ge 1}{R_2}^j z^j\Big) \\
&=&
\sum_{n\ge 1} \sum_{\substack{i+j = n\\i,j\ge 1}}
\big(R_1^i z^i\big)\cdot \big(R_2^j z^j\big)\\
&=&
\sum_{n\ge 1} \Big(\sum_{\substack{i+j = n\\i,j\ge 1}}
{R_1}^i {R_2}^j\Big)   z^n\\
 &\unlhd&
\sum_{n\ge 1}(R_1+R_2)^n z^n
\  =\  \bA_{R_1+R_2}(z).
\end{eqnarray*}
For composition, 
letting $R_0= 1+R_1+R_2$:
\begin{eqnarray*}
\big(\bA_{R_1}\circ\bA_{R_2}\big)(z)
&\unlhd& \big(\bA_{{R_0}}\circ\bA_{{R_0}}\big)(z)
\ =\  \sum_{i\ge 1} {R_0}^i\Big( \sum_{j\ge 1} {R_0}^j z^j\Big)^i\\
&=& \sum_{i\ge 1} \Big( \sum_{j\ge 1} {R_0}^{1+j} z^j\Big)^i
\ \unlhd\  \sum_{i\ge 1}\Big(\sum_{j\ge 1}({R_0}^2 z)^j\Big)^i\\
&\unlhd& \sum_{n\ge 1}(2{R_0}^2 z)^n \ =\ \bA_{2 {R_0}^2}(z).
\end{eqnarray*}
\end{proof}

\subsection{Bounded operators}

The main tool for showing that the solution $w=\bT$ to $w=\Theta(w)$ has
a positive radius of convergence, which is essential to employing the
methods of analysis, is to show that $\Theta$ is 
{\em bounded}.

\begin{definition}
For $R>0$ define the simple operator $\bA_R$ by
\[
\bA_R(w)\ =\ \sum_{j\ge 1}R^j\,w^j.
\]
An operator $\Theta$ is \underline{bounded} if  
$\big(\exists R>0\big)\Big(\Theta(w)\ole\bA_R(z + w)\Big)$, that is,
\[
\big(\exists R>0\big)\,\big(\forall \bT\in\dom\big)\,\Big(\Theta(\bT)\unlhd \bA_R(z + \bT)\Big).
\]
\end{definition}

Of course we will want to use integer values of $R$ when working with integral
operators.

\subsection{When is an elementary operator bounded?} 

The properties {\em weakly retro} and {\em monotone} investigated 
earlier hold for all elementary operators. This is certainly not the 
case with the bounded property. In this subsection we give a 
simple univariate test for being bounded. 

As mentioned before, any
$\bE\in\dom[z,w]$ plays a dual role in this paper, one as a bivariate
power series and the other as an elementary operator. Each of these 
roles has its own definition as to what bounded means, namely:
\begin{eqnarray*}
\bE\unlhd \bA_R(z+w)&\Leftrightarrow&
\big(\forall i,j\ge 1\big)\,\Big(e_{ij} \le [z^i\,w^j]\,\bA_R(z+w)\Big)\\
\bE\ole\bA_R(z+w)&\Leftrightarrow&
\big(\forall \bT\in\dom[z]\big)\,\Big(\bE(z,\bT)\,\unlhd\,\bA_R(z+\bT)\Big).
\end{eqnarray*}
The two definitions are equivalent.

\begin{lemma}\label{two defs}
Let $\bE$ be an elementary operator.
\begin{thlist}
\item
$\bE(z,w)$ is bounded as an operator iff\, $\bE(z,z)$ is bounded
as a power series.
Indeed
\begin{eqnarray*}
\bE(z,w)\ole\bA_R(z+w)
&\Rightarrow& 
\bE(z,z)\unlhd\bA_{2R}(z)\quad\text{for }R>0\\
\bE(z,z)\unlhd\bA_R(z)&\Rightarrow& 
\bE(z,w)\ole\bA_{R}(z+w) \quad\text{for }R>1.
\end{eqnarray*}
\item
The equivalence of bivariate bounded and operator bounded
follows from
\begin{eqnarray*}
\bE(z,w)\unlhd\bA_R(z+w)&\Rightarrow& \bE(z,w)\ole\bA_R(z+w)\quad\text{for }R>0\\
\bE(z,w)\ole\bA_R(z+w)&\Rightarrow& \bE(z,w)\unlhd\bA_{2R}(z+w)\quad\text{for }R>1.
\end{eqnarray*}
\end{thlist}
\end{lemma}
\begin{proof}
For (a) suppose $R>0$ and $\bE(z,w)\ole\bA_R(z+w)$.
Since $z \in\dom[z]$, we have 
\[
\bE(z,z) \ \unlhd\  \sum_{j\ge 1} R^j (2z)^j \ =\ \bA_{2R}(z),
\]
so $\bE(z,z)$ is a bounded power series.

Conversely, suppose $R>1$ and $\bE(z,z)\unlhd \bA_R(z)$.
Then
\[
\bE(z,z)\ \unlhd\ \sum_{j\ge 1} R^j z^j,
\]
so for $n\ge 1$
\[
[z^j]\, \bE(z,z) \ \le \ R^j.
\]
Then from $\bE(z,w) =\sum e_{i,j} z^i w^j$
we have $e_{i,j} \le R^{i+j}$,
so
\begin{eqnarray*}
\bE(z,w)&\unlhd &\sum_{i,j\ge 1} R^{i+j} z^i w^j\\
&\unlhd& \sum_{i,j\ge 1} R^{i+j} \binom{i+j}{i} z^iw^j\\
& =& \bA_R(z+w).
\end{eqnarray*}
Applying Lemma \ref{bidom} gives
$\bE(z,w)\ole\bA_R(z+w)$. 

For (b) the first claim is just Lemma \ref{bidom}. For the second
claim suppose $R>1$ and 
$\bE(z,w)\ole\bA_R(z+w)$.
From the first part of (a) we have 
$\bE(z,z)\unlhd \bA_{2R}(z)$ and then from the second part 
$\bE(z,w)\ole\bA_{2R}(z+w)$.

\end{proof}

\begin{corollary} \label{cs bdd}
Given $\bA\in\dom[z]$, the constant operator $\Theta_\bA$ 
as well as the simple operator $\bA(w)$ are bounded iff $\rho_\bA >0$.
\end{corollary}

\subsection{Bounded operators form a closed set} \label{bdd closed sect}

\begin{lemma} \label{bdd is closed}
The set of bounded operators is  closed.
\end{lemma}
\begin{proof}
Let $\Theta,\Theta_1,\Theta_2$ be bounded operators as witnessed by
the following:\\
$\Theta(w)\ole \bA_R(z+w)$, $\Theta_1(w)\ole \bA_{R_1}(z+w)$ 
and $\Theta_2(w)\ole \bA_{R_2}(z+w)$.
With $c>0$ we have from Lemma \ref{bounded updir}
\begin{eqnarray*}
\big(c \Theta\big)(w)& \ole& c \bA_R(z+w)\ \ole\ \bA_{(1+c)R}(z+w)\\ 
\big(\Theta_1 + \Theta_2\big)(w)& \ole& \bA_{R_1}(z+w) + \bA_{R_2}(z+w)\ \ole\ \bA_{R_1+R_2}(z+w)\\
\big(\Theta_1  \Theta_2\big)(w)
& \ole& \bA_{R_1}(z+w)  \bA_{R_2}(z+w)\ \ole\ 
\bA_{R_1+R_2}(z+w)\\
\big(\Theta_1 \circ\Theta_2\big)(w)& \ole& \bA_{R_1}(z+w) \circ\bA_{R_2}(z+w)\ \ole\ \bA_{2(1+R_1+R_2)^2}(z+w).
\end{eqnarray*}
\end{proof}

\begin{lemma} All restrictions of standard operators are 
bounded operators.
\label{std bdd}
\end{lemma}
\begin{proof}
Let $\Delta$ be a standard operator. Then for any $\bbM\subseteq \bbP$ we
have $\Delta_\bbM \ole \Delta$, so it suffices to show the standard
operators are bounded. But this is evident from the well known fact that
\begin{eqnarray*}
\MSet(w)
&\ole& \Cycle(w)\  
\ole\ \DCycle(w)\\ 
&&\ole\ 
\Seq(w)\ =\ \sum_{n\ge 1}w^n\ = \bA_1(w)\ \ole\ \bA_1(z+w).
\end{eqnarray*}
So the choice of $R$ is $R=1$.
\end{proof}


\subsection[I]{When dominance of operators gives dominance of fixpoints} 

This is part of proving that the solution $w=\bT$ to $w=\Theta(w)$
has a positive radius of convergence.

\begin{lemma}\label{domi}
Let $\bT_i$ satisfy the recursion equation
$\bT_i\, =\,\Theta_i(\bT_i)$ for $i=1,2$\,. 
If the $\Theta_i$ are retro operators,
$\Theta_1\,\ole\,\Theta_2$, and
$\Theta_1$ or $\Theta_2$ is monotone 
then $\bT_1\,\unlhd\,\bT_2$\,.
\end{lemma}

\begin{proof}
Since each $\Theta_i(w)$ is a retro operator, 
by Lemma \ref{retro lem} we have 
\[
\bT_i\ =\ \lim_{n\rightarrow\infty}{\Theta_i}^n(0)\,.
\]
Let us use induction to show
\[
{\Theta_1}^n(0)\,\unlhd\,{\Theta_2}^n(0)
\]
holds for $n\ge 1$\,.
For $n=1$ this follows from the assumption that $\Theta_2$
dominates $\Theta_1$. So suppose it holds for $n$. Then
\[
\begin{array}{l}
{\Theta_1}^{n+1}(0)\ 
\unlhd\ \Theta_1\Big({\Theta_2}^n(0)\Big)\ 
\unlhd\ {\Theta_2}^{n+1}(0)\quad\text{if $\Theta_1$ is monotone}\\
{\Theta_1}^{n+1}(0)\ 
\unlhd\ \Theta_2\Big({\Theta_1}^n(0)\Big)\ 
\unlhd\ {\Theta_2}^{n+1}(0)\quad\text{if $\Theta_2$ is monotone}.
\end{array}
\]
Thus $\bT_1\,\unlhd\,\bT_2$\,.
\end{proof}

\subsection{The nonzero radius lemma}
To apply complex analysis methods to a solution $\bT$ of a recursion
equation we need $\bT$ to be analytic at 0. 

\begin{lemma}\label{modify bound}
Let $\Theta$ be a retro operator with 
$\Theta(w)\sqsubseteq\bA_R(z+w)$. 
Then
\[
\Theta(w)\ \sqsubseteq\ \bA_R(z+w) - Rw.
\]
\end{lemma}
\begin{proof}
Since $\Theta$ is retro there is a sequence $\sigma_n$ of functions such
that for $\bT\in\dom[z]$,
\[
[z^n]\,\Theta(\bT)\ =\ \sigma_n\big(t_1,\ldots,t_{n-1}\big). 
\]
Let 
\[
\Phi(w) \ :=\ \sum_{n\ge 2}R^n (z+w)^n, 
\]
which is easily seen to be a retro operator. 
Choose $\widehat{\sigma}_n$ such that
for $\bT\in \dom[z]$ 
\[
[z^n]\,\Phi(\bT)
\ =\ 
\widehat{\sigma}_n\big(t_1,\ldots,t_{n-1}\big).
\]
Then, since $\bA_R(z+w) = R(z+w) + \Phi(w)$, from the dominance of 
$\Theta(w)$ by $\bA_R(z+w)$ we have, for any $t_i\ge 0$
and $n\ge 2$,
\[
\sigma_n\big(t_1,\ldots,t_{n-1}\big)\ \le\ R t_n\, + \,
\widehat{\sigma}_n\big(t_1,\ldots,t_{n-1}\big).
\]
As the left side does not depend on $t_n$ we can put $t_n=0$ 
to deduce
\[
\sigma_n\big(t_1,\ldots,t_{n-1}\big)\ \le\ 
\widehat{\sigma}_n\big(t_1,\ldots,t_{n-1}\big),
\]
which gives the desired conclusion.
\end{proof}

\begin{lemma} \label{pos rad lemma}
Let $\Theta$ be a bounded retro operator. Then $w=\Theta(w)$ has
a unique solution $w=\bT$, and $\rho_\bT>0$.
\end{lemma}
\begin{proof}
By Lemma \ref{retro lem} we know there is a unique solution $\bT$. Choose
$R>1$ such that $\Theta(w) \sqsubseteq \bA_R(z + w)$. 
From Lemma \ref{modify bound} we can change this to
\begin{equation} \label{theta R}
\Theta(w)\ \sqsubseteq\ \bA_R(z + w) - Rw. 
\end{equation}
The right side is a monotone retro operator, 
so Lemma \ref{domi} says that the
fixpoint $\bS$ of $\bA_R(z+w)-Rw$ dominates the fixpoint $\bT$ of $\Theta(w)$. 
Let
\[
\bS\ =\ \bA_R(z+\bS) - R\bS.
\]
To show $\rho_\bT>0$
it suffices to show $\rho_\bS>0$.
We would like to sum the geometric series $\bA_R\big(z+\bS(z)\big)$; 
however since we do not yet know that $\bS$ is analytic at $z=0$ we 
perform an equivalent maneuver by multiplying both sides of equation
\eqref{theta R}
by $1-Rz - R\bS$ to obtain the quadratic equation
\[
 \left( R+{R}^{2} \right) {\bS}^{2}+ \left( {R}^{2}z+Rz-1 \right) \bS+Rz
\ =\ 0.
\]
The discriminant of this equation is
\[
\bD(z)\ =\ \left( {R}^{2}z+Rz-1 \right)^2\, -\, 
4\left( R+{R}^{2} \right) Rz.
\]
Since $\bD(0) = 1 $ is positive 
it follows that $\sqrt{\bD(z)}$ is analytic in a 
neighborhood of $z=0$. Consequently $\bS(z)$ has a nonzero radius
of convergence.

\end{proof}

\subsection{The set of composite operators}

The sets $\cO_E$ and $\cO_I$ of operators that we 
eventually will exhibit as 
``guaranteed to give the universal law'' will be subsets of the 
following {\em composite} operators.

\begin{definition} \label{composite}
The composite operators are those obtained from 
the base operators, namely
\begin{thlist}\itemsep=1ex
\item
the elementary operators $\bE(z,w)$ and 
\item
the $\bbM$-restrictions of the standard operators:
$\MSet_\bbM$, 
$\Cycle_\bbM$, 
$\DCycle_\bbM$ 
and $\Seq_\bbM$,
\end{thlist}
using the variables $z,w$, scalar multiplication by positive reals,
 and the binary operations addition 
\mbox{\rm($+$)}, multiplication \mbox{\rm($\cdot$)} 
and composition \mbox{\rm($\circ$)}.
\end{definition}

\begin{lemma} \label{2 constr props}
The set of composite operators is closed under the arithmetical operations
and
all composite operators $\Theta$ are monotone and weakly retro.
\end{lemma}
\begin{proof}
The closure property is immediate from the definition of the set of
composite operators,
the monotone property is from Lemma \ref{monotone closed}, and
the weakly retro property is from Lemma \ref{retro is full} (b).
\end{proof}

\begin{figure}[h]
\centerline{\psfig{figure=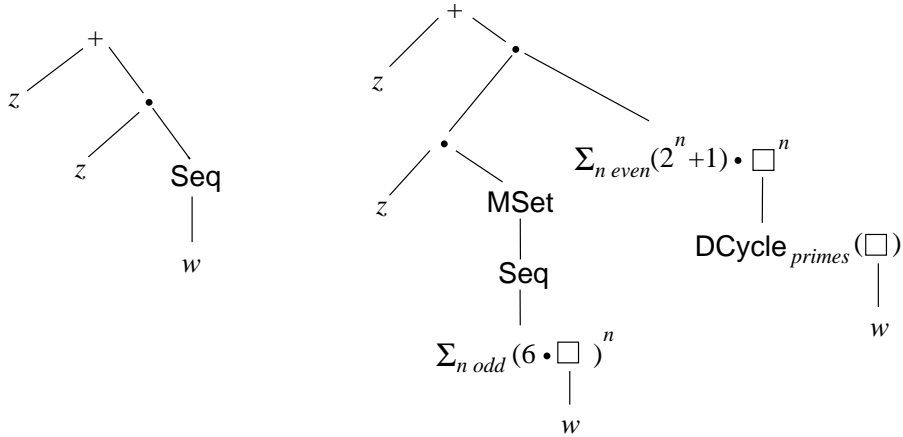}}
\caption{Two examples of term trees \label{term tree}}
\end{figure}

An expression like $z\, +\,z \Seq(w)$ that describes
how a composite operator is constructed is 
called a \emph{term}. Terms can be visualized as trees, for example
the term just described and the term in \eqref{ex1} have the 
trees shown in Figure \ref{term tree}. (A small empty box in the figure
 shows where the argument below the box is to be inserted.)
Composite operators are, like their counterparts called {\em term functions}
in universal algebra and logic, valued for the fact that one has the possibility
to (1) define functions on the class by induction on terms, and (2) 
one can prove facts about the class by induction on terms.

Perhaps the simplest explanation of why we like the composite operators 
$\Theta$ so much is: we have a routine procedure to convert the
equation $w = \Theta(w)$ into an equation $w = \bE(z,w)$ where $\bE$
is elementary. This is the next topic.

\subsection{Representing a composite operator $\Theta$ at $\bT$}

In order to apply analysis to the solution $w=\bT$ of a recursion
equation $w=\Theta(w)$ we want to put the
equation into the form $w = \bE(z,w)$ with $\bE$ analytic on $\bT$. 
The next definition describes a natural candidate for
$\bE$ in the case that $\Theta$ is composite.

\begin{definition} \label{repr}
Given a base operator $\Theta$ and a $\bT\in\dom[z]$ 
define an elementary operator 
$\bE^{\Theta,\bT}$ 
as follows:
\begin{thlist}
\item
$\bE^{\bE,\bT} \,=\, \bE$ for $\bE$ an elementary operator.
\item
For $\Theta=\MSet_\bbM$ let
$\bE^{\Theta,\bT}\,=\,\sum_{m\in\bbM} \bZ\big(\sS_m,w,\bT(z^2),\ldots,\bT(z^m)\big)$. 
\item
For $\Theta=\DCycle_\bbM$ let
$\bE^{\Theta,\bT}\,=\,\sum_{m\in\bbM} \bZ\big(\sC_m,w,\bT(z^2),\ldots,\bT(z^m)\big)$. 
\item
For $\Theta=\Cycle_\bbM$ let
$\bE^{\Theta,\bT}\,=\,\sum_{m\in\bbM} \bZ\big(\sD_m,w,\bT(z^2),\ldots,\bT(z^m)\big)$. 
\item
For $\Theta=\Seq_\bbM$ let
$\bE^{\Theta,\bT}\,=\,\sum_{m\in\bbM} w^m$.
\end{thlist}

Extend this to all composite operators using the obvious inductive definition:
\begin{eqnarray*}
\bE^{c \Theta,\bT}&:=&c \bE^{\Theta,\bT}\\
\bE^{\Theta_1 + \Theta_2,\bT}& :=& \bE^{\Theta_1,\bT} + \bE^{\Theta_2,\bT}\\
\bE^{\Theta_1  \Theta_2,\bT}& :=& \bE^{\Theta_1,\bT}  \bE^{\Theta_2,\bT}\\
\bE^{\Theta_1 \,\circ\, \Theta_2,\bT}& :=& 
\bE^{\Theta_1,\Theta_2(\bT)}\big(z,\bE^{\Theta_2,\bT}\big).
\end{eqnarray*}
\end{definition}
The definition is somewhat redundant as 
the $\Seq_\bbM$ operators are included in the elementary operators.

\begin{lemma}\label{repres}
For $\Theta$ a composite operator and $\bT\in\dom[z]$ we have
\[
\Theta(\bT)\ =\ \bE^{\Theta,\bT}(z,\bT).
\]
We will simply say that $\bE^{\Theta,\bT}$ 
\underline{represents $\Theta$ at $\bT$}. 
\end{lemma}
\begin{proof}
By induction on terms.
\end{proof}

\subsection{Defining linearity for composite operators} 

\begin{definition}
Let $\Theta$ be a composite operator. We say $\Theta$ is 
\underline{linear} (in $w$) if the elementary operator 
$\bE^{\Theta,z}$ representing $\Theta$ at $z$ is
linear in $w$. Otherwise we say $\Theta$ is \underline{nonlinear} (in $w$).
\end{definition}

\begin{lemma} \label{lin inv}
Let $\Theta$ be a composite operator. Then the elementary operator 
$\bE^{\Theta,\bT}(z,w)$ representing $\Theta$ at $\bT$ is either 
linear in $w$ for all $\bT\in\dom[z]$, or it is 
nonlinear in $w$ for all $\bT\in\dom[z]$.
\end{lemma}
\begin{proof}
Use induction on terms.
\end{proof}

\subsection{When $\bT$ belongs to $\dom^\star[z]$}

\begin{proposition}\label{dom star}
Let $\Theta$ be a bounded nonlinear retro composite operator. Then there
is a unique solution $w = \bT$ to $w = \Theta(w)$, and $\bT\in\dom^\star[z]$,
that is, $\rho_\bT\in(0,\infty)$ and $\bT(\rho_\bT)<\infty$.
\end{proposition}
\begin{proof}
From Lemma \ref{pos rad lemma} we know that $w=\Theta(w)$ has a unique
solution $\bT\in\dom[z]$, and $\rho :=\rho_\bT>0$.
Let $\bE(z,w)$ be the elementary operator representing $\Theta$ at $\bT$.
Then $\bT = \bE(z,\bT)$. As $\Theta$ is nonlinear there is a positive
coefficient $e_{ij}$ of $\bE$ with $j\ge 2$. Clearly
\[
\bT(x) \ \ge\ e_{ij} x^i \bT(x)^j\quad\text{for }x\ge 0.
\]
Divide through by $\bT(x)^2$ and take the limsup of both sides as $x$
approaches $\rho^-$ to see that $\bT(\rho) < \infty$, and thus $\rho<\infty$.
This shows $\bT\in\dom^\star[z]$.
\end{proof}

\subsection{Composite operators that are open for $\bT$}

Many examples of elementary operators enjoy the open property, but 
(restrictions of) the standard operators rarely do: only the various
$\Seq_\bbM$ and $\Delta_{\{1\}}$ for $\Delta$ any of the standard
operators. 

For the standard operators other than $\Seq$, and hence for
most of the composite operators, it is 
very important that we use the concept of `open at $\bT$'
when setting up for the Weierstra{\ss} Preparation Theorem.

\begin{definition} \label{op open def}
Let $\bT\in\dom^\star[z]$.
A composite operator $\Theta$ is \underline{open for $\bT$} iff
$\bE^{\Theta,\bT}$ is 
open at $\big(\rho,\bT(\rho)\big)$.
\end{definition}

The next lemma determines when the base operators are open for a given 
$\bT\in\dom^\star[z]$.

\begin{lemma} \label{open constr}
Suppose $\bT\in\dom^\star[z]$ and let $\rho\in(0,\infty)$ be its 
radius of convergence.
Then the following hold:
\begin{thlist}
\item
An elementary operator $\bE$ is open for $\bT$ iff it is open
at $\big(\rho,\bT(\rho)\big)$.
\item
A constant operator $\Theta_\bA(w)$ is open for $\bT$\, iff\,
$\rho<\rho_\bA$.
\item
A simple operator $\bA(w)$ is open for $\bT$ \,iff\, 
$\bT(\rho)<\rho_\bA$. 
\item
$\Seq_\bbM$ is open for $\bT$ iff $\bbM$ is finite or
$\bT(\rho)<1$.
\item
$\MSet_\bbM$ is open for $\bT$\, iff\, 
$\bbM=\{1\}$ or $\rho<1$.
\item
$\DCycle_\bbM$, or $\Cycle_\bbM$, is open for $\bT$\, iff\, 
$\bbM=\{1\}$\ or\  
$(\bbM$ is finite and $\rho<1)$ or 
$(\bbM$ is infinite and $\rho,\bT(\rho)<1)$. 
\end{thlist}
\end{lemma}

\begin{proof}
For (a) note that an open operator represents itself at $\bT$.
For (b) and (c) use Lemma \ref{cs open}.
For (d) note that $\Seq_\bbM(w)$
is the simple operator $\bA(w) :=\sum_{m\in\bbM} w^m$, so (c) applies.

For (e) let $\bE := \bE^{\Theta,\bT}$
where $\Theta := \MSet_\bbM$.
Then 
\[
\bE(z,w)\ :=\ \sum_{m\in\bbM} 
\bZ\big(\sS_m,w,\bT(z^2),\ldots,\bT(z^m)\big).
\]
If $\bbM=\{1\}$ 
then $\bE(z,w)=w$ and (c) applies.
So suppose $\bbM\neq\{1\}$.
The term $\bT(z^2)$ appears in $\bE(z,w)$, and this diverges
at $\rho + \varepsilon$ if $\rho \geq 1$. Thus $\rho < 1$ is 
a necessary condition for $\bE$ to be open for $\bT$. 

So suppose $\rho<1$. 
The representative for $\MSet$ 
dominates the representative of any
$\MSet_\bbM$. Thus for any $x\in(0,\sqrt{\rho})$ and $y>0$:
\[
\bE(x,y)\ \le\ e^y \exp\Big(\sum_{m\ge 2}\bT\big(x^m\big)\big/m\Big)\ <\ \infty.
\]
Since one can find $\varepsilon>0$ such that the 
right hand side is finite at
$\big(\rho+\varepsilon,\bT(\rho)+\varepsilon\big)$, it follows
that $\MSet_\bbM$ is open for $\bT$ when $\rho<1$.
\medskip

For (f) let
$\bE := \bE^{\Theta,\bT}$
where $\Theta := \DCycle_\bbM$.
Then 
\begin{eqnarray*}
\bE(z,w)
& :=& 
\sum_{m\in \bbM} \bZ\big(\sC_m,w,\bT(z^2),\ldots,\bT(z^m)\big)\\
& =& 
\underbrace{\sum_{m\in\bbM}\frac{1}{m} w^m}_{\bA(w)}\ +\ 
\underbrace{\sum_{k\ge 2}\frac{\varphi(k)}{k}\sum_{jk\in\bbM}\frac{1}{j}\bT(z^k)^j}_{\bB(z)}.
\end{eqnarray*}
If $\bbM=\{1\}$ then, as before, there are no further restrictions needed
as $\bE(z,w) := w$. 
So now suppose $\bbM\neq\{1\}$.
The presence of some $\bT(z^k)$ with $k\ge 2$ in the expression 
for $\bE(z,w)$ shows, as in (e), that a necessary condition is $\rho<1$.
This condition implies $\rho_\bB \ge \sqrt{\rho}$.

If $\bbM$ is finite then $\rho_\bA=\infty$, and $\rho_\bB \ge \sqrt{\rho}$, 
consequently $\bE$ is open at $\big(\rho,\bT(\rho)\big)$.

If $\bbM$ is infinite then $\rho_\bA=1$.
Suppose $\bE$ is open at $\big(\rho,\bT(\rho)\big)$.
Then $\bA\big(\bT(\rho)+\varepsilon\big)$ converges for some
$\varepsilon>0$, so $\bT(\rho)<1$. The conditions $\rho,\bT(\rho)<1$
are easily seen to be sufficient in this case.

For the $\Cycle_\bbM$ case let
$\bE := \bE^{\Theta,\bT}$
where $\Theta := \Cycle_\bbM$.
\begin{eqnarray*}
\Cycle_\bbM\big(\bT(z)\big)& =&  \frac{1}{2}\,\DCycle_\bbM\big(\bT(z)\big)\\
&& +\ 
\frac{1}{4}
\sum_{m\in\bbM} 
\begin{cases}
2\bT(z) \bT(z^2)^{(m-1)/2} &\text{if $m$ is odd}\\
\bT(z)^2 \bT(z^2)^{(m-2)/2}\ +\ \bT(z^2)^{m/2} &\text{if $m$ is even}.
\end{cases}
\end{eqnarray*}
Thus
\begin{eqnarray*}
\bE(z,w)
& :=& 
\sum_{m\in \bbM} \bZ\big(\sD_m,w,\bT(z^2),\ldots,\bT(z^m)\big)\\
& =& 
\frac{1}{2}\sum_{m\in\bbM}\frac{1}{m} w^m\ +\ 
\frac{1}{2}\sum_{k\ge 2}\frac{\varphi(k)}{k}\sum_{jk\in\bbM}\frac{1}{j}\bT(z^k)^j\\
&& +\ 
\frac{1}{4}\sum_{m\in\bbM} 
\begin{cases}
2 w \bT(z^2)^{(m-2)/2} &\text{if $m$ is odd}\\
w^2 \bT(z^2)^{(m-1)/2}\ +\ \bT(z^2)^{m/2} &\text{if $m$ is even}
\end{cases}
\end{eqnarray*}
and we can use the same arguments as for $\DCycle$\,.
\end{proof}

\subsection{Closure of the composite operators that are open for 
$\bT\in\dom^\star[z]$} \label{aeo}

\begin{lemma} \label{aeo lem}
Suppose $\bT\in\dom^\star[z]$.  Then the following hold:
\begin{thlist}
\item
The set of composite operators that are open for $\bT$
is closed under addition, scalar multiplication and multiplication.
\item
Given composite operators $\Theta_1,\Theta_2$ with 
$\Theta_2$ open for $\bT$ and $\Theta_1$ open for
$\bT_1 := \Theta_2(\bT)$, the composition $\Theta_1\circ \Theta_2$ 
is open for $\bT$.
\end{thlist}
\end{lemma}

\begin{proof}
Just apply Lemma \ref{closure of open E}.
\end{proof}

\subsection{Closure of the composite integral operators that are open 
for $\bT\in\idom^\star[z]$} \label{aio}

\begin{definition}
  $\idom^\star[z] = \idom[z] \cap \dom^\star[z]$.
\end{definition}

\begin{lemma} \label{aio lem}
Suppose $\bT\in\idom^\star[z]$.  Then the following hold:
\begin{thlist}
\item
The set of integral composite operators that are open for $\bT$
is closed under addition, positive integer scalar multiplication and multiplication.
\item
Given integral composite operators $\Theta_1,\Theta_2$ with 
$\Theta_2$ open for $\bT$ and $\Theta_1$ open for
$\bT_1 := \Theta_2(\bT)$, the composition $\Theta_1\circ \Theta_2$ 
is integral and open for $\bT$.
\end{thlist}
\end{lemma}

\begin{proof}
This is just a repeat of the previous proof, noting that at
each stage we are dealing with integral operators acting on $\idom[z]$.
\end{proof}

\subsection{A special set of operators called $\cO$}

This is the penultimate step in describing the promised collection 
of recursion equations.

\begin{definition} \label{def cO}
Let $\cO$ be the set of operators that can be constructed from
\begin{thlist}\itemsep=1ex
\item
the bounded and open elementary operators $\bE(z,w)$ and 
\item
the $\bbM$-restrictions of the standard operators:
$\MSet_\bbM$, 
$\Cycle_\bbM$, 
$\DCycle_\bbM$ 
and $\Seq_\bbM$,
where in the case of the cycle constructions we require the set $\bbM$ to
be either finite or to satisfy $\sum_{m\in\bbM} 1/m \,=\, \infty$,
\end{thlist}
using the variables $z,w$, scalar multiplication by positive reals,
 and the binary operations addition 
\mbox{\rm($+$)}, multiplication \mbox{\rm($\cdot$)} 
and composition \mbox{\rm($\circ$)}.

Within $\cO$ let $\cO_E$ be the set of bounded and open 
elementary operators; and let $\cO_I$ be the closure 
under the arithmetical operations of the
bounded and open integral elementary operators along with the
standard operators listed in (b). 
\end{definition}

Clearly $\cO$ is a subset of the composite operators.

\begin{lemma} \label{O facts}
\begin{thlist}
\item
Every $\Theta\in\cO$ is a bounded monotone and weakly retro operator.
\item
Each of the sets $\cO,\cO_E,\cO_I$ is closed under the 
arithmetical operations.
\end{thlist}
\end{lemma}
\begin{proof}
For (a) 
we know from our assumption on the elementary operators in $\cO$ and
Lemma \ref{std bdd} that the base operators in $\cO$ are bounded---then 
Lemma \ref{bdd is closed} shows that all members of $\cO$ are 
bounded. 
All members of $\cO$ are monotone and weakly retro
by Lemma \ref{2 constr props}. 
Regarding (b), use Lemma \ref{closure of open E} (b) for $\cO_E$, and 
Definition \ref{def cO} for the other two sets.
\end{proof}

\begin{lemma}\label{open integral}
Let $\Theta\in\cO_I$. 
If $\,\bT\in\idom^*[z]$ and 
$\Theta(\bT)(\rho_\bT)<\infty$ then $\Theta$ is open for $\bT$.
\end{lemma}
\begin{proof}
Since $\bT\in\idom^\star$ we must have $\rho:=\rho_\bT < 1$.
Let 
\[
\cO^\star\ :=\ \{\Theta\in\cO_I : \Theta(\bT)(\rho) < \infty\}.
\]
An induction proof will show that for $\Theta\in\cO^\star$
we have $\Theta$ open for $\bT$.
The elementary base operators of $\cO^\star$ are given to be
open, hence they are open for $\bT$.
The restrictions of the standard operators in $\cO^\star$
are covered by parts (d)--(f) of Lemma \ref{open constr}, with
one exception.
We need to verify in certain $\DCycle$
and $\Cyc$ cases that $\bT(\rho)<1$. In these cases
one has $\bbM$ infinite, and then we must have $\bT(\rho)<1$
in order for $\Theta(\bT)$ to converge at $z=\rho$ since $\sum_{m\in\bbM} 1/m \,=\, \infty$.

For the induction step simply apply Lemma \ref{aio lem}.
\end{proof}

\subsection{The Main Theorem}

The following is our main theorem, exhibiting many $\Theta$ for which
$w = \Theta(w)$ is a recursion equation whose solution satisfies 
the universal law.  Several examples follow the proof. 

\begin{theorem}\label{Main Thm}
Let $\Theta_1$ be a nonlinear retro member of 
$\cO_E$, respectively $\cO_I$, and let $\bA(z)\in \dom[z]$,
respectively $\bA(z)\in\idom[z]$, be such that 
$\bA(\rho_\bA)=\infty$.
Then there is a unique $\bT\in\dom[z]$, respectively $\bT\in\idom[z]$,
 such that $\bT = \bA(z) + \Theta_1(\bT)$.
The coefficients of $\bT$ satisfy the
universal law $\pmb{(\star)}$ in the form
\[
t(n)\ \sim\ q\sqrt
{\frac{\rho\bE_z\big(\rho,\bT(\rho)\big)} 
{2\pi \bE_{ww}\big(\rho,\bT(\rho)\big)}} 
\cdot
\rho^{-n} n^{-3/2}
\quad\text{for } n\equiv d \mod q. 
\]
Otherwise $t(n)=0$.
Thus $\pmb{(\star)}$ holds on $\{n : t(n) > 0\}$.
The constants $d,q$ are from the
shift periodic form $\bT(z) = z^d\bV(z^q)$.
\end{theorem}

\begin{proof}
Let $\Theta(w) = \bA(z) + \Theta_1(w)$, by Lemma \ref{O facts} 
a member of $\cO_E$, respectively
$\cO_I$.
By Proposition \ref{dom star}
there is a unique solution $w=\bT$ to $w=\bA(z) + \Theta(w)$ and 
$\bT\in\dom^\star[z]$. 
Let $\bE_1(z,w)=\bE^{\Theta,\bT}$. Then the elementary representative $\bE$
of $\bA(z)+\Theta(w)$ is given by
\[
\bE(z,w) \ :=\ \bA(z) + \bE_1(z,w).
\]
We will verify the hypotheses (a)--(e) of Theorem \ref{basic thm}.

$\bT = \bE(z,\bT)$ by Lemma \ref{repres};
this is \ref{basic thm} (a).
The fact that $\bT\in\dom^\star[z]$ is \ref{basic thm} (b).
By Lemma \ref{lin inv} we get \ref{basic thm} (c).
Since $\bA(0) = 0$ and $\bA \neq 0$
it follows that $\bA_z \neq 0$. 
As $\bE(z,0) = \bA(z)$ it follows that $\bE_z \neq 0$. 
This is \ref{basic thm} (d).

To show $\Theta$ is open for $\bT$ we note that in the
case of the operators coming from $\cO_E$ they are given 
to be open elementary operators; and for the case they
are coming from $\cO_I$ use Lemma \ref{open integral}.
This gives \ref{basic thm} (e).
\end{proof}

\subsection{Applications of the main theorem} \label{applications}

One readily checks that all the recursion equations given in
Table \ref{std examp} satisfy the hypotheses of Theorem \ref{Main Thm}.
One can easily produce more complicated examples such as
\[
w\ =\ 3z^3\,+\,z^4 \Cycle(w)
\,+\,w^2 \DCycle(w) \,+\,\MSet_2(w).
\]

Such simple cases barely scratch the surface 
of the possible applications of Theorem \ref{Main Thm}. Let us turn to 
the more dramatic example given early in \eqref{ex1}, namely:
\[
w\ =\  z \ +\ z \MSet
\Big(\Seq\big(
\sum_{n\in\Odd} 6^n w^n
\big)\Big)  
\sum_{n\in\Even}(2^n+1)
\big(\DCycle_{\Primes}(w)\big)^n\,.
\]
We will analyze this from `the inside out', naming the operators encountered
as we work up the term tree. First we give names to the nodes of the term tree:
\[
\begin{array}{l @{\quad}l@{\quad} l@{\quad} l } 
\Phi_1\ :=\ \sum_{n\in\Odd} 6^n w^n&
\Phi_2\  :=\  \DCycle_{\Primes}(w)&
\Phi_3\  :=\  \Seq\big(\Phi_1\big)\\
\Phi_4\  :=\  \MSet\big(\Phi_3\big)&
\Phi_5\  :=\  \sum_{n\in\Even}(2^n+1) w^n&
\Phi_6\  :=\  \Phi_5(\Phi_4)\\
\bA(z)\ :=\ z&
\Theta_1\  :=\  z\Phi_4 \Phi_6&&.
\end{array}
\]
Now we argue that each of these operators is in $\cO_I$:
\begin{thlist}
\item
$\Phi_1$ is an elementary (actually simple) integral operator with radius 
of convergence $1/6$. Thus it is bounded.  Since it diverges at its radius 
of convergence, it is open. Thus $\Phi_1\in\cO_I$.

\item
$\Phi_2$ is a restriction of $\DCycle$ to the 
set of prime numbers; since $\sum_{m\in\Primes} 1/m = \infty$ we have
$\Phi_2\in\cO_I$.

\item
$\Phi_3$ is in $\cO_I$ as it is a composition of two operators in $\cO_I$.

\item
$\Phi_4$ is in $\cO_I$ as it is a composition of two operators in $\cO_I$.

\item
$\Phi_5$ is an elementary (actually simple) integral operator with radius 
of convergence $1/2$. Thus it is bounded.  Since it diverges at its radius 
of convergence, it is open. Thus $\Phi_5\in\cO_I$.

\item
$\Phi_6$ is in $\cO_I$ as it is a composition of two operators in $\cO_I$.

\item
$\Theta_1$ is in $\cO_I$ as it is a product of two operators in $\cO_I$.

\item
$\Theta_1$ is a nonlinear retro operator in $\cO_I$.
\end{thlist}

Thus we have an equation 
$w\ =\ \bA(z) + \Theta_1(w)$
that satisfies the hypotheses of Theorem \ref{Main Thm}; consequently
the solution $w=\bT(z)$ has coefficients satisfying the universal law.

\subsection{Recursion specifications for planar trees}

When working with either labelled trees or planar trees the
recursion equations are elementary. Here is a popular
example that we will examine in detail.

\begin{example}[Planar Binary Trees] The defining equation is 
\[
w \ =\  z \,+\, z w^2.
\]
This simple equation can be handled directly since it is a quadratic,
giving the solution
\[
\bT(z) \ =\ \frac{1 - \sqrt{1-4z^2}}{2z}.
\]
Clearly $\rho = 1/2$ and for $n\ge 1$ we have $t(2n)=0$, and Lemma
\ref{binomial} gives
\begin{eqnarray*}
t(2n-1) 
&=& (-1)^n \frac{4^n}{2}\binom{1/2}{n}
 \sim\ 
 \frac{4^n}{2}
\cdot \frac{n^{-3/2}}{2 \sqrt{\pi}}
\ =\ \frac{1}{\sqrt{\pi}} 4^{n-1} n^{-3/2}.
\end{eqnarray*}

For illustrative purposes let us examine this in light of 
the results in this paper.  Note that 
\[
\bE(z,w)\ :=\  z + z w^2
\]
is in the desired form $\bA(z) + \Theta_1(w)$ with 
$\bA(\rho_\bA) = \infty$
and
$\Theta_1$ a bounded retro nonlinear (elementary) operator.

The constants $d,q$ of the shift periodic form are given by: 
\begin{thlist}
\item
$d=1$ as $\bE_0(z) = z$ implies $E_0=\{1\}$.
\item
$q = 2$ as $E_0 = \{1\}$, $E_2 = \{1\}$, and otherwise $E_j=\text{\O}$;
thus $\bigcup E_n+(n-1)d = (E_0 - 1)\cup(E_2+1) = \{0,2\}$, so 
$q = \gcd\{0,2\}= 2$. 
\end{thlist}
Thus $t(n)>0$ implies $n\equiv 1 \mod 2$, that is, $n$ is an odd number.
For the constant in the asymptotics we have 
\begin{eqnarray*}
\bE_z(z,w)&=& 1+w^2\\
\bE_{ww}(z,w)&=& 2z.
\end{eqnarray*}
In this case we know $\rho = 1/2$ and $\bT(\rho)=1$ (from solving the 
quadratic equation), so
\begin{eqnarray*}
\bE_z\big(\rho,\bT(\rho)\big)&=& 2\\
\bE_{ww}\big(\rho,\bT(\rho)\big)&=& 1.
\end{eqnarray*}
Thus
\begin{eqnarray*}
t(n)&\sim&
q 
\sqrt{\frac{\rho\bE_z\big(\rho,\bT(\rho)\big)} 
{2\pi \bE_{ww}\big(\rho,\bT(\rho)\big)}} 
\cdot \rho^{-n} n^{-3/2}\\
&=& 2\sqrt{\frac{1}{2\pi}} \cdot 2^{n} n^{-3/2}\\
&=& \sqrt{\frac{2}{\pi}} \cdot 2^n n^{-3/2}\quad\text{for }
n\equiv 1 \mod 2.
\end{eqnarray*}
\end{example}

\subsection{On the need for integral operators}
Since the standard operators, and their restrictions, are defined
on $\dom[z]$ it would be most welcome if one could unify the treatment
so that the main theorem was simply a theorem about operators on $\dom[z]$
instead of having one part for elementary operators on $\dom[z]$, and another
part for integral operators acting on $\idom[z]$. However the following 
example indicates that one has to exercise some caution when working
with standard operators that mention $\bT(x^j)$ for some $j\ge 2$.

Let 
\[
\Theta(w)\ :=\ \frac{z}{2}\big(1 + \MSet_2(w)\big).
\]
This is 1/2 the operator one uses to define (0,2)-trees. This operator
is clearly in $\cO$ and of the form 
$\bA(z) + \Theta_1(w)$; 
however it is not in either $\cO_E$ or $\cO_I$, as required by the main
theorem.

$\Theta$ is clearly retro and monotone.
Usual arguments show that $w = \Theta(w)$ has a unique solution $w=\bT$
which is in $\dom^\star$, and we have
\begin{equation}\label{one}
\bT(\rho) 
\  =\  \frac{1}{2}\rho + \frac{1}{4}\rho\bT(\rho)^2 + 
\frac{1}{4}\rho\bT(\rho^2).
\end{equation}
Since $\Theta(\bT)$ involves $\bT(z^2)$ it follows that
$\rho\le 1$ (for otherwise $\bT(\rho^2)$ diverges).

Suppose $\rho<1$.
Following P\'olya let us write the equation for $\bT$ as $w = \bE(z,w)$
where 
\[
\bE(z,w)\ :=\  
 \frac{1}{2}z + \frac{1}{4}z w^2 + \frac{1}{4}z \bT(z^2).
\]
Then the usual condition for the singularity $\rho$ is 
$1 = \bE_w(\rho,\bT(\rho))$, that is
\begin{equation}\label{two}
1\  =\ \frac{1}{4}(\rho 2\bT(\rho))\ =\ \frac{\rho\bT(\rho)}{2},
\end{equation}
so $\rho\bT(\rho)=2$.

Putting $\bT(\rho) = 2/\rho$ into equation \eqref{one} gives
\[
\frac{2}{\rho}\  =\  \frac{1}{2}\rho + \frac{1}{\rho} +
\frac{1}{4}\rho \bT(\rho^2),
\]
so
\[
4\  =\  2\rho^2 +  \rho^2\bT(\rho^2). 
\]
Since  $\bT(\rho^2) < \bT(\rho)=2/\rho$ we have
\[
4\ <\ 2\rho^2 + 2\rho,
\]
a contradiction as $\rho<1$.

Thus $\rho = 1$, and we cannot apply the method of P\'olya since
$\bE(z,w)$ is not holomorphic at $\big(1,\bT(1)\big)$.

\section{Algorithmic Aspects}

\subsection{An algorithm for nonlinear}

Given a term $\Phi(z,w)$ that describes a composite operator $\Theta$
there is a simple algorithm to determine if $\Theta$ is nonlinear. 
Let us use the abbreviation $\Delta$ for the various standard unary 
operators and their restrictions as well as the elementary operators 
$\bE(z,w)$\,.
We can assume that any 
occurrence of a $\bbM$-restriction of a standard operator $\Delta$ in $\Phi$ 
is such that 
$\bbM\,\neq\,\{1\}$ since if 
$\bbM\,=\,\{1\}$ then $\Delta_\bbM$ is just the identity operator. 
Given an occurrence of a $\Delta$ in $\Phi$ let $T_\Delta$ be the full
subtree of $\Phi$ rooted at the occurrence of $\Delta$\,.

\noindent
{\sf An algorithm to determine if a composite $\Theta$ is nonlinear}
\begin{itemize}
\item
First we can assume that constant operators $\Theta_\bA$ are only located at the
leaves of the tree.
\item
If there exists a $\Delta$ in the tree of $\Phi$ such that a leaf $w$ is 
below $\Delta$\,,  
where $\Delta$ is either a restriction of a standard
operator or a nonlinear elementary $\bB$, then $\Theta$ is nonlinear.
\item
If there exists a node labelled with multiplication  in the tree such that 
each of the two branching nodes have a $w$ on or below them then 
$\Theta$ is nonlinear.
\item
Otherwise $\Theta$ is linear in $w$.
\end{itemize}
\begin{proof}[Proof of the correctness of the Algorithm.]
(A routine induction argument on terms.)
\end{proof}

\section{Equations $w = \bG\big(z,w\big)$ with mixed 
sign coefficients}\label{mixed}

\subsection{Problems with mixed sign coefficients} \label{sub mixed}

We would like to include the possibility of mixed sign
coefficients in a recursion equation $w = \bG(z,w)$. 
The following table shows the key steps we used to prove 
$\pmb{(\star)}$ holds in the nonnegative case, and the situation
if we try the same steps in the mixed sign case.
\[
\begin{tabular}{|l |l |l|}
\hline
&$\bG\in\bR^{\ge 0}[[z,w]]$& $\bG\in\bR[[z,w]]$\\
&Nonnegative $\bG$&Mixed Signs $\bG$\\
\hline
Property&Reason &Reason\\
\hline
$(\exists!\bT)\,(\bT=\bG(z,\bT)$&$g_{01}=0$& $g_{01}=0$\\ 
$\rho>0$&$\bG$ is bounded& $\bG$ is abs. bounded\\
$\rho<\infty$&$\bG$ is nonlinear in $w$& (?)\\
$\bT(\rho)<\infty$&$\bG$ is nonlinear in $w$&(?)\\
$\bG$ holomorphic in nbhd of $\bT$& $\bG(\rho+\varepsilon,\bT(\rho)+\varepsilon)<\infty$&(?)\\
$\bG_{ww}\big(\rho,\bT(\rho)\big)\neq 0$&$\bG$ is nonlinear in $w$& (?)\\
$\bG_z\big(\rho,\bT(\rho)\big)\neq 0$&$\bG_0(z)\neq 0$& (?)\\
$\DSing = \{z : z^q = \rho^q\}$&$\Spec\,\bG_w\big(z,\bT(z)\big)$ is nice& (?)\\
\hline
\end{tabular}
\]

As indicated in this table, many of the techniques that we used for the 
case of a nonnegative equation do not carry over to the mixed case. 
\begin{thlist}
\item
To show that a unique solution $w=\bT$ exists in the mixed sign case
we can use the retro property, precisely as with the nonnegative case. 
The condition for $\bG\in \bbR[[z,w]]$ to be retro is that $g_{01}=0$. 
\item
To show $\rho>0$ in the nonnegative case we used the existence of an 
$R>0$ such that $\bE(z,\bT) \unlhd \bA_R\big(z+\bT\big)$. In the mixed 
sign case we could require that $\bG(z,\bT)$ be absolutely dominated 
by $\bA_R\big(z+\bT\big)$.
\item
To show $\rho < \infty$ and $\bT(\rho)<\infty$ we used the nonlinearity 
of $\bE(z,w)$ in $w$. Then $\bT = \bE(z,\bT)$ implies 
$\bT(x) \ge e_{ij}x^i\bT(x)^j$ for some $e_{ij}>0$ with $j\ge 2$. 
This conclusion does not follow in the mixed sign case.
\item
After proving that $\bT\in\dom^\star[z]$,
to be able to invoke the theoretical machinery of $\S$\,\ref{theor sect}
we required that $\bE$ be open at $\big(\rho,\bT(\rho)\big)$, that is,
\[
\big(\exists \varepsilon>0\big)\,
\Big(\bE\big(\rho+\varepsilon,\bT(\rho)+\varepsilon\big) \,<\, \infty \Big).
\]
This shows $\bE$ is holomorphic on a neighborhood of $\bT$.
In the mixed sign case there seems to be no such easy condition unless we
know that $\sum_{ij} |g_{ij}|\rho^i\bT(\rho)^j < \infty$.
\item
In the nonnegative case, if $\bE$ is nonlinear in $w$ then $\bE_{ww}$ does
not vanish, and hence it cannot be 0 when evaluated at 
$\big(\rho,\bT(\rho)\big)$. With the mixed signs case, proving 
$\bG_{ww}\big(\rho,\bT(\rho)\big)\neq 0$
requires a fresh analysis.
\item
A similar discussion applies to showing $\bG_z\big(\rho,\bT(\rho)\big)\neq 0$.
\item
Finally there is the issue of locating the dominant singularities.
The one condition we have to work with is that the dominant singularities
must satisfy $\bG_w\big(z,\bT(z)\big) = 1$.
In the nonnegative case we were able to use the analysis of the
spectrum of $\bE_w$:
\[
\Spec\Big(\bE_w\big(z,\bT(z)\big)\Big)\ = \ \bigcup_n E_n + (n-1)\odot T.
\]
This tied in with an expression for the spectrum of $\bE\big(z,\bT(z)\big)$.
However for the mixed case we only have
\[
\Spec\Big(\bE_w\big(z,\bT(z)\big)\Big)\ \subseteq \ \bigcup_n E_n + (n-1)\odot T.
\]

In certain mixed sign equations one has a promising property, namely
\[
\bG_w\big(z,\bT(z)\big)\in \dom[z].
\]
This happens with the equation for identity trees.
In such a case put  
$\bG_w\big(z,\bT(z)\big)$ in its pure periodic form $\bU(z^p)$. Then 
the necessary condition on the dominant singularities $z$ becomes 
simply
\[
z^p \ =\ \rho^p.
\]
If one can prove $p=q$, as we did with elementary recursions, 
then the dominant singularities are as simple as one could hope for.
\end{thlist}
There is clearly considerable work to be done to develop a theory of solutions
to mixed sign recursion equations.  

\subsection{The operator $\Set$} \label{Set sect} \label{sub set}
The above considerations led us to omit the 
popular $\Set$ operator from our list of standard
combinatorial operators.
In the equation $w = z+z\Set(w)$ for the class of {\em identity trees}
(that is, trees with only one automorphism), one can readily show that the only
dominant singularity of the solution $\bT$ is $\rho$. But if we look at more complex
equations, like 
\[
w = z +z^3 + z^5 + z \Set\big(\Set(w)\MSet(w)\big),
\]
the difficulties of determining the locations of the dominant singularities 
appear substantial.

\begin{example}

Consider the restrictions of the $\Set$ operator 
\begin{eqnarray*}
\Set_\bbM(\bT)& =& \sum_{m\in\bbM}\Set_m(\bT),\quad\text{where}\\
\Set_m(\bT)& = & \bZ\big(\sS_m,\bT(z),-\bT(z^2),\ldots,(-1)^{m+1}\bT(z^m)\big). 
\end{eqnarray*}
Thus in particular
\[
\Set_2(\bT)\ = \ \frac{1}{2}\big(\bT(z)^2 - \bT(z^2)\big).
\]

The recursion equation
\[
w\ =\ z\,+\,\Set_2(w)
\]
exhibits different behavior than what has been seen so far 
since the solution is
$
\bT(z)\, =\, z,
$
which is not a proper infinite series. 
The solution certainly does not have coefficients 
satisfying the universal law, nor does it have a finite radius of convergence
that played such an important role. 

We can modify this equation slightly to obtain a more interesting solution, 
namely let
\[
\Theta(w)\ =\ z\,+\,z^2\,+\,\Set_2(w).
\]
Then $\Theta$ is integral retro, and the unique solution $w=\bT$ to 
$w = \Theta(w)$ is
\[
\bT(z)\ =\ z \,+\, z^2\, +\,z^3\,+\,z^4\,+\,2z^5\,+\,3z^6\,+\,6z^7\,+\,11z^8\,+\,\cdots
\]
with $t(n)\ge 1$ for $n\ge 1$.
Consequently we have the radius of convergence 
\begin{equation}\label{rho in 01}
\rho :=\rho_\bT\in[0,1].
\end{equation}
We will give a detailed proof that $\bT$ has coefficients satisfying the 
universal law, to hint at the added difficulties that might occur
in trying to add $\Set$ to our standard operators. 

Let 
\[
\Theta_1(w) \ =\  z + z^2 + \frac{1}{2}w^2
\]
a bounded open nonlinear retro elementary operator, hence an 
operator in $\cO_E$ to which the Main Theorem applies. 
For $\bA\in\idom$ note that
$\Theta(\bA)\, \unlhd\,\Theta_1(\bA)$, so we can use the monotonicity of $\Theta_1$
to argue that $\Theta^n(0)\,\unlhd\,{\Theta_1}^n(0)$ for all $n\ge 1$. 
Thus $\bT$ is dominated by the solution
$\bS$ to $w=\Theta_1(w)$. At this point we know that 
$\rho_\bT \ge\rho_\bS> 0$. 

Since $\rho\in(0,1]$ 
we have $\bT(x^2) < \bT(x)$ for $x\in(0,\rho)$. Thus for $x\in(0,\rho)$
\[
\bT(x)\ >\ x + x^2 + \frac{1}{2}\big(\bT(x)^2 - \bT(x)\big),
\]
or
\[
\frac{3}{2}\bT(x)\ >\ x + x^2 + \frac{1}{2}\bT(x)^2. 
\]
Thus $\bT(x)$ cannot approach $\infty$ as $x\rightarrow \rho^-$. 
Consequently $\bT(\rho)<\infty$, and then we must also have $\rho<1$.
By defining
\begin{eqnarray*}
\bG(z,w)&:=& z\,+\,z^2+\, \frac{1}{2}\big(w^2 - \bT(z^2)\big)\\
& =& \frac{1}{2}w^2\,+\,z\,+\,\frac{1}{2}z^2
\,-\,\frac{1}{2}z^4
\,-\,\frac{1}{2}z^6 \,-\,\cdots,
\end{eqnarray*}
we have the recursion equation $w = \bG(z,w)$ satisfied by $w=\bT$, and
$\bG(z,w)$ has mixed signs of coefficients.

As $\rho<1$ we know that $ \bG(z,w)$
is holomorphic in a neighborhood of the graph of $\bT$, so
a necessary condition for $z$ to be a dominant singularity is that
$\bG_w\big(z,\bT(z)\big)=1$, that is
$\bT(z)\ =\ 1$.
Since $\bT$ is aperiodic,
this tells us we have a {\em unique dominant singularity}, namely $z=\rho$, 
and we have $\bT(\rho)\ =\ 1$.

Differentiating the equation
\[
\bT(z)\ =\ z + z^2 + \frac{1}{2}\bT(z)^2 - \frac{1}{2}\bT(z^2)
\]
gives
\[
\bT'(z)\ =\ 1 + 2z + \bT(z)\bT'(z) - z \bT'(z^2)
\]
or equivalently
\[
\big(1-\bT(z)\big)\bT'(z)
\ =\ 
(1 + 2z)\,-\, z\bT'(z^2)
\quad\text{for }|z|<\rho.
\]
Since $\rho<1$ we know that 
\[
\lim_{\substack{z\rightarrow \rho\\|z|<\rho}} \Big((1+2z) -z\bT'(z^2)\Big)
\ =\ (1 + 2\rho) - \rho\bT'(\rho^2).
\]
Let $\lambda$ be this limiting value. 
Consequently
\begin{equation}\label{TT'}
\lim_{\substack{z\rightarrow \rho\\|z|<\rho}} \big(1-\bT(z)\big)\bT'(z)\ =\ \lambda.
\end{equation}
By considering the limit along the real
axis, as $x\rightarrow \rho^-$, we see that $\lambda\ge 0$, so
\[
(1 + 2\rho) - \rho\bT'(\rho^2)\ =\ \lambda \ \ge \ 0.
\]
Let 
\[
\bF(z,w)\ :=\ w \,-\,\Big(z + z^2 + \frac{1}{2}w^2 - \frac{1}{2}\bT(z^2)\Big). 
\]
Then 
\[
\bF_z(z,w)\ =\ -\Big(1 + 2z - z\bT'(z^2)\Big)\ =\ z\bT'(z^2) - (1+2z),
\]
so
\[
\bF_z\big(\rho,\bT(\rho)\big)\ =\  \rho\bT'(\rho^2) - (1+2\rho)\ =\ -\lambda.
\]
If $\lambda > 0$ then 
$\bF_z\big(\rho,\bT(\rho)\big)\, <\, 0$; and since $\bF_{ww} = -2$ we have
\[
\bF_z\big(\rho,\bT(\rho)\big)
\bF_{ww}\big(\rho,\bT(\rho)\big)\ > \ 0.
\]
This means we have all the hypotheses needed to apply Proposition \ref{crucial1}
to get the square root asymptotics which lead to the universal law for $\bT$.

To conclude that we indeed have the universal law we will show that $\lambda >0$.
Let $\alpha\in[\rho,1]$. Then for $x\in(0,\rho)$ we have 
$\bT(x^2) \le \alpha \bT(x)$, and thus for $x \in (0,\rho)$
\[
\bT(x)\ >\ x + x^2 + \frac{1}{2}\Big(\bT(x)^2 - \alpha\bT(x)\Big).
\]
Let 
\[
\bU(x)\ =\ x + x^2 + \frac{1}{2}\Big(\bU(x)^2 - \alpha\bU(x)\Big).
\]
Then 
\begin{eqnarray*}
\bU(x)& =& \frac{1}{2} \Big( (2+\alpha) - \sqrt{ (2+\alpha)^2 - 8(x + x^2)}\Big)\\
\rho_\bU &=& -\frac{1}{2} + \frac{1}{4}\sqrt{4 + 2(2+\alpha)^2}.
\end{eqnarray*}

Now for $x\in I := (0,\min(\rho,\rho_\bU))$
\[
(2+\alpha)\bT(x) - \bT(x)^2
\ >\ 
(2+\alpha)\bU(x) - \bU(x)^2
\]
so
\[
\bU(x)^2 - \bT(x)^2\ >\ (2 + \alpha)\big( \bU(x) - \bT(x)\big).
\]
Thus $\bU(x) \neq \bT(x)$ for $x\in I$. If $\bU(x)>\bT(x)$ on $I$ then
\[
\bU(x) + \bT(x)\ >\ 2 + \alpha \quad\text{for }x\in I.
\]
But this is impossible since on $I$ we have 
\begin{eqnarray*}
\bU(x) &<& \bU(\rho_\bU) = 1 + \alpha/2\\ 
\bT(x) &<& \bT(\rho)\ =\ 1.
\end{eqnarray*}
Thus we have
\[
\bU(x) \ <\ \bT(x)\quad\text{on } I.
\]
If $\rho_\bU \le \rho$ then $\bU(\rho_U) \le \bT(\rho_U)$, which is also impossible. 
Thus
\[
\rho < \rho_\bU.
\]

Now define a function $f$ on $[\rho,1]$ that maps $\alpha\in[\rho,1]$ to 
$\rho_\bU$ as given in the preceding lines, that is:
\[
f(\alpha)\ =\ (-1/2) + (1/4)\sqrt{4 + 2(2+\alpha)^2}.
\]
Then $\alpha \in [\rho,1]$ implies $f(\alpha) \in (\rho,1]$. 
Calculation gives $f^3(1) = 0.536\ldots$, so 
\[
\rho \ <\ 0.54
\]
Since $\rho_\bS = \big(\sqrt{3}-1\big)/2 = 0.366\ldots$ 
we have $\rho^2 < \rho_\bS$, 
and then
\[
\bT'(\rho^2)\ <\ \bS'(\rho^2),
\]
so
\[
-\lambda \ =\ \rho\bT'(\rho^2) - (1+2\rho)\ <\ \rho\bS'(\rho^2) - (1+2\rho)\ <\ 0
\]
since $x\bS'(x^2) - (1+2x)\ < \ 0$ for $x\in (0,0.55)$.
This proves $\lambda>0$, and hence the universal law holds for the coefficients of
$\bT$.

This example shows that the generating function $\bT^*$ for the class
of identity (0,1,2)-trees satisfies the universal law. We have
$\bT^*$ defined by the equation $w = z + z * \Set_{\{1,2\}}(w)$, and it turns out  
that $t_n^* = t_{n+1}$. (We discovered this connection with $\bT^*$
when looking for the first few coefficients of $\bT$ in the {\em On-Line 
Encyclopedia of Integer Sequences}.)
\end{example}

\begin{example}
In the 20 Steps paper of 
Harary, Robinson and Schwenk \cite{HRS} 
the asymptotics for the class of {\em identity trees} (those with no nontrivial automorphism)
was successfully analyzed
by first showing that the associated recursion equation
\[
w\ =\ z \,+\, z \Set(w) 
\]
has a unique solution $w=\bT\in\dom^\star[z]$. 
Then
\[
\bG(z,w)\ :=\ z\,+\,z e^w\cdot\exp\Big(\sum_{m\ge 2}(-1)^{m+1}\bT(z^m)/m\Big)
\]
is holomorphic in a neighborhood of the graph of $\bT$. One has 
\[
z + \bG_w(z,w)\ =\ \bG(z,w),
\]
so the necessary condition 
$\bG_w\big(z,\bT(z)\big)=1$ 
for a dominant singularity 
is just the condition $\bT(z) = 1+z$.
$\rho$ is the only solution of this equation on the circle 
of convergence as $\bT-z\,\unrhd\,0$ and is aperiodic. 
Consequently {\em the only dominant singularity is $z=\rho$}.

The equation for identity trees is of mixed signs.
From
\begin{equation}\label{SetEq}
\bT(z)\ =\ z \prod_{j\ge 1}\big(1+z^j\big)^{t_j}
\end{equation}
we can calculate the first few values of $t(n)$ for identity trees:\footnote{One
can also look up sequence number A004111 in the 
{\em On-Line Encyclopedia of Integer Sequences.}}
\[
\begin{array}{| c | c | c | c | c | c | c | c | c | c |}
t(1)&t(2)&t(3)&t(4)&t(5)&t(6)&t(7)&t(8)&t(9)&t(10)\\
\hline
1&1&1&2&3&6&12&25&52&113
\end{array}
\]
Returning to the definition of $\bG$ we have
\begin{eqnarray*}
\bG(z,w)& =& z + z e^w
\sum_{n\ge 0}
\Big(\sum_{m\ge 2} (-1)^{m+1}\big(t(1)z^m + t(2)z^{2m} +\cdots)\big/m\Big)^n\Big/ n!\\
& =& z + z e^w
\sum_{n\ge 0}\Big(-\big(z^2 + z^{4} +\cdots)\big/2
\ +\ \big(z^3 + z^{6} +\cdots)\big/3\\
&&\qquad\ -\ \big(z^4 + z^{8} +\cdots)\big/4
\Big)^n\Big/ n!\\
& =& z + z e^w
\sum_{n\ge 0}\Big(-z^2/2 + z^3/3 -3 z^{4}/4 +\cdots\Big)^n\Big/ n!\\
& =&z + z  e^w
\Big(1 -z^2/2 + z^3/3 - 5z^{4}/8 +\cdots\Big).
\end{eqnarray*}
Thus for some of the $z^iw^j$ the coefficients are positive, and some are
negative; $\bG(z,w)$ is a mixed sign operator.
\end{example}

If one were to form more complex operators $\Theta$ by adding the operator 
$\Set$ and its restrictions $\Set_\bbM$ to our set of Standard Operators, then 
there is some hope for proving that one always has the universal
law holding for the solution to $w = \Theta(w)$, provided one has a solution
that is not a polynomial. 

The hope stems from the fact that although the $\bG(z,w)$
associated with $\Theta$ may have mixed sign coefficients, when it comes to the
condition $\bG_w\big(z,\bT(z)\big) = 1$ on the dominant singularities we have 
the good fortune that $\bG_w\big(z,\bT(z)\big)\unrhd 0$, that is, it expands 
into a series with nonnegative coefficients. The reason is quite simple, 
namely using the bivariate generating function we have
\[
\Set_m(\bT)\ =\ [u^m]\,\exp\Big(\sum_{n\ge 1} (-1)^{n-1}u^n\bT(x^n)/n\Big).
\]
Letting $\bG_m(z,w)$ be $\bZ_m\big(\sS_m,w,-\bT(z^2),\ldots,(-1)^{m+1}\bT(z^m)\big)$
we have
\[
\bG_m(z,w)\ =\ 
[u^m]\,e^{uw}\cdot\exp\Big(\sum_{n\ge 2} (-1)^{n-1}u^n\bT(x^n)/n\Big),
\]
thus
\begin{eqnarray*}
\frac{\partial \bG_m}{\partial w}
& =& \frac{\partial }{\partial w}
[u^m]\,e^{uw}\cdot\exp\Big(\sum_{n\ge 2} (-1)^{n-1}u^n\bT(x^n)/n\Big)\\
& =& [u^m]\,u e^{uw}\cdot\exp\Big(\sum_{n\ge 2} (-1)^{n-1}u^n\bT(x^n)/n\Big)\\
& =& [u^{m-1}]\, e^{uw}\cdot\exp\Big(\sum_{n\ge 2} (-1)^{n-1}u^n\bT(x^n)/n\Big)\\
&=& \bG_{m-1}(z,w).
\end{eqnarray*}

Consequently if we put $\bG_w\big(z,\bT(z)\big)$ into 
its pure periodic form $\bU(z^p)$ then we have the necessary condition 
$z^p = \rho^p$ on the dominant singularities $z$. Letting 
$z^d\bV(z^q)$ be the shift periodic form of $\bT(z)$ it follows that 
$q \big| p$. If we can show that $p = q$ then 
$\DSing = \{z : z^q = \rho^q\}$, which is as simple as possible.
Indeed this has been the case with the few examples we have worked out by 
hand.


\section{Comments on Background Literature} \label{history}

Two important sources offer global views on finding asymptotics. 

\subsection{The ``20 Step algorithm'' of \cite{HRS}\label{20 step}} 

This 1975 paper by Harary, Robinson and Schwenk 
is in good part a heuristic for
how to apply P\'olya's method\footnote{The paper also has a 
proof that the generating function for the class of identity trees 
(defined by $w = z + z \Set(w)$) satisfies the universal law.}, 
and in places the 
explanations show an affinity for operators close to the
original ones studied by P\'olya. 
For example it says that 
$\bG(z,w)$ should be analytic for $|z|\,<\,\sqrt{\rho_\bT}$ and
$|w|\,<\,\infty$. This strong condition on $w$ fails for most of the simple
classes studied by Meir and Moon, and hence for the setting of this
paper.

The algorithm of {\em 20 Steps} also discusses how to find asymptotics for the class of 
free trees obtained from a rooted class defined by recursion. 
Given a class $\cT$ of rooted trees let $\cU$ be the associated class of 
free (unrooted) trees, that is, the members of $\cU$ are the same 
as the members of
$\cT$ except that the designation of an element as the root has been
removed. Let the corresponding generating series\footnote{{\em 20 Steps}
uses $t(n)$ to denote the number of rooted trees in $\cT$ on $n+1$ points, whereas
$u(n)$ denotes the number of free trees in $\cU$ on $n$ points. We will let $t(n)$
denote the number of free trees in $\cT$ on $n$ points, as we have 
done before, since this will have no material effect on the efficacy of 
the 20 Steps.} 
be $\bT(z)$ and $\bU(z)$. 

The initial assumptions are only two: that $\bT$ is not a polynomial and 
it is aperiodic. 
Step 2 of the 20 steps is: express $\bU(z)$ in terms of $\bT(z)$ and $\bT(z^2)$.
{\em 20 Steps}
says that Otter's {\em dissimilarity characteristic} can usually be 
applied to achieve this.
Step 20 is to deduce that $u(n) \sim C \rho^{-n} n^{-5/2}$.

This outline suggests that it is widely possible to find the asymptotics
of the coefficients $u(n)$, and evidently this gives a second universal law
involving the exponent $-5/2$ instead of the $-3/2$. 
Our investigations suggest that determining the growth 
rate of the associated classes of free trees will be quite challenging.

Suppose $\cT$ is a class of rooted trees for which the P\'olya 
style analysis has been successful, that we have found the radius of convergence
$\rho \in (0,1)$, that $\bT(\rho)< \infty$, and that 
$t(n) \sim C \rho^{-n}  n^{-3/2}$. What can we say about the generating
function $\bU(z)$ for the corresponding class of free trees?

Since we have the inequality $t(n)/n \le u(n) \le t(n)$ (note
that one has only $n$ ways to choose the root in a free tree of size $n$), it follows
by the Cauchy-Hadamard Theorem that $\bU(z)$ has the same radius of convergence 
$\rho$ as $\bT(z)$. From this and the asymptotics for $t(n)$ it also follows 
that one can find $C_1,C_2>0$ such that for $n\ge 1$
\[
C_1 \rho^{-n}  n^{-5/2}\ \le \ u(n)\ \le \  C_2 \rho^{-n}  n^{-3/2}.
\]
Thus $u(n)$ is sandwiched between a $-5/2$ expression and a $-3/2$ expression.
In the case that $\cT$ is the class of {\em all} rooted trees, Otter \cite{Otter}
showed that
\[
\bU(z)\ =\ z\Big(\bT(z) - \Set_2(\bT)\Big)
\ =\  z\Big(\bT(z) - \frac{1}{2}\big(\bT(z)^2 - \bT(z^2)\big)\Big),
\]
and from this he was able to find the asymptotics for $u(n)$ with 
a $-5/2$ exponent.

However let $\cT$ be the class of rooted
trees such that every node has either 2 or 5 descending branches. The 
recursion equation for $\bT(z)$ is
\[
\bT(z)\ =\ z\,+\,z\Big(\MSet_2(\bT) \,+\, \MSet_5(\bT)\Big),
\]
and $\bT(z)$ is aperiodic.
By Theorem \ref{Main Thm} we know that the coefficients of $\bT(z)$ satisfy the
universal law
\[
t(n)\ \sim\ C \rho^{-n}  n^{-3/2}.
\]

Let $\cU$ be the corresponding set of free trees. Note that when one converts a
rooted tree $T$ in $\cT$ to a free tree $F$, a root with 2 descending branches 
will give a node of degree 2 in $F$, and a root with 5 descending branches will 
give a node of degree 5 in $F$. 
Any non root node with 2 descending branches will give a node of degree 3 in $F$;
and any non root node with 5 descending branches will give a node of degree 6 in $F$.
Thus $F$ will have exactly one node of degree 2 or degree 3, and not both, so one
can identify the node that corresponds to the root of $T$. This means that there
is a bijection between the rooted trees on $n$ vertices in $\cT$ and the free
trees on $n$ vertices in $\cU$. Consequently $t(n)=u(n)$, and thus
\[
u(n)\ \sim\ C \rho^{-n}  n^{-3/2}.
\]
Clearly $u(n)$ cannot also satisfy a $-5/2$ law. Such examples are easy
to produce. 

Thus it is not clear to what extent the program of {\em 20 Steps} can be carried
through for free trees. It seems that free trees are rarely defined by a 
single recursion equation, and it is doubtful if there is always a 
recursive relationship between $\bU(z)$ and $\bT(z),\bT(z^2),\ldots$. 
Furthermore it is not clear what the possible asymptotics for the $u(n)$ 
could look like---is it possible that one will always have either a 
$-3/2$ or a $-5/2$ law? Since a class $\cU$ derived from a $\cT$ which 
has a nice recursive specification can be defined by
a monadic second order sentence, there is hope that the $u(n)$ 
will obey a reasonable asymptotic law. (See Q5 in $\S\,$\ref{Open Probs}.)

In {\em 20 Steps} consideration is also given to techniques for calculating 
the radius 
of convergence $\rho$ of $\bT$ and the constant $C$ that appears in the
asymptotic formula for the $t(n)$\,. In this regard the reader should consult
the paper of Plotkin and Rosenthal \cite{Pl:Ro} as there are evidently some
numerical errors in the constants calculated in {\em 20 Steps}.

\subsection{Meir and Moon's global approach}\label{MM sec}

In 1978 Meir and Moon \cite{Me:Mo1}
considered classes $\cT$ of trees with generating 
functions $\bT(z)\,=\,\sum_{n\ge 1} t(n) z^n$
such that
\begin{quote} 
\begin{thlist}\itemsep=1ex
\item[1]
$t(1)\,=\,1$\,;
\item[2]
$\cT$ can be obtained by taking certain forests of trees
from $\cT$ and adding a root to each one
(this choice of certain forests is evidently a `construction');
\item[3]
this `construction' and `conditions implicit in the definition' of $\cT$
give rise to a `recurrence relation' for the $t(n)$\,, evidently 
a sequence $\sigma$ of functions $\sigma_n$ such that
$t(n)\,=\,\sigma_n\big(t(1),\ldots,t(n-1)\big)$\,;
\item[4]
there is an `operator' $\Gamma$\,, acting on (possibly infinite) sequences
of power series, such that the recurrence relation for the $t(n)$ 
`can be expressed in terms of generating series', for example
$\bT(z)\,=\,\Gamma\big(\bT(z),\bT(z^2),\ldots\big)$\,, which is abbreviated
to $\bT(z)\,=\,\Gamma\{\bT(z)\}$\,.
\end{thlist}
\end{quote}
This is the most penetrating presentation we have seen of a foundation 
for recursively defined classes of trees, a goal that we find most
fascinating since to prove global results one needs a global setting. 
However their conditions have limitations that we want to point out.

(1) requires $t(1)=1$, so there is only one object of size 1; this means   
multicolored trees are ruled out. (2) indicates that one
is using a specification\footnote{We have not needed a specification
language so far in this paper---for this comment it is useful. 
Let $\bullet$ denote the tree with one node, and $\bullet\big/\Box$ says to
add a root to any forest in $\Box$\,.
}
like
$\cT\,=\,\{\bullet\}\,\cup\,\bullet\big/\Seq_\bbM(\cT)$\,. The recurrence
relation in (3) is the one item that seems to be appropriately general. It 
corresponds to what we call `retro'. (4) is too
vague; after all, a function of $\big(\bT(z),\bT(z^2),\ldots\big)$ 
is really just a function of $\bT(z)$ since $\bT(z)$ completely 
determines all the $\bT(z^k)$. This
formulation is surely motivated by the desire to include the $\MSet$ 
construction; perhaps the authors were thinking of `natural' functions
of these arguments like $\sum_{n\ge 1} \bT(z^n)/n$\,; or perhaps something 
of an effective nature, an algorithm.

After this general discussion, without any attempt to prove theorems in this
context, they turn the focus to \emph{simple} classes $\cT$ of rooted trees,
namely simple classes are those for which
\begin{quote}
\begin{itemize}
\item[(M1)]
the generating series $\bT(z)$ is defined by a `simple' recursion 
equation
\[
w\ =\  z \bA(w)\,,
\]
where $\bA\,\in\,\bbR^{\ge 0}[[z]]$ with $\bA(0)\,=\,1$\,. 
\end{itemize}
\end{quote}
Additional conditions\footnote{Their original 1978 conditions 
had a minor restriction, that $a(1) > 0$.  That was soon replaced 
by the condition (M3)---see for example \cite{Me:Mo3}.
}
are needed to prove their theorems, 
namely
\begin{quote}
\begin{itemize}
\item
[(M2)] $\bA(w)$ is analytic at 0,
\item
[(M3)] $\gcd\{n\in \bbP: a(n)>0\}\, =\,1$\,, 
\item
[(M4)] $\bA(w)$ is not a linear polynomial $aw+b$\,, and
\item
[(M5)] $\bA(w)\,=\,w \bA'(w)$ has a positive solution
$y\,<\,\rho_\bA$
\end{itemize}
\end{quote}
to guarantee that the methods of P\'olya apply to give the 
asymptotic form $\pmb{(\star)}$\,. 
(M2) makes $\bT(z)$ analytic at 0, (M3) ensures that $\bT(z)$ is aperiodic,
(M4) leads to $\rho_\bT\,< \,\infty$ and $\bT(\rho_\bT)\,<\,\infty$\,,
and (M5) shows $\bF(z,w)\,:=\,w\,-\,z \bA(w)\,=\,0$
is holomorphic in a neighborhood
of $\big(\rho_\bT,\bT(\rho_\bT)\big)$.

Thanks to the restriction to recursion identities based on simple operators
they are able to employ the more powerful condition (M5) instead of our 
condition $\bA(\rho_\bA)\,=\,\infty$\,. Our condition 
is easier to use in practice, and it covers the two examples frequently cited
by Meir and Moon, namely planar trees with $\bA(w)\,=\,\sum_{n\ge 0}w^n$
and planar binary trees with $\bA(w)\,=\,1+w^2$\,.
For the simple recursion equations one can replace 
our $\bA(\rho_\bA)\,=\,\infty$ by the condition
$\bA'(\rho_\bA)\,=\,\infty$.

\section{Open Problems} \label{Open Probs}
\begin{thlist}
\item [Q1]
If $\bT = \bE(z,\bT)$ with $\bE\in\dom[z,w]$ and $\bT(z)$ has 
the shift periodic form $z^d \bV(z^q)$ then one can use the 
spectrum calculus to show there is an $\bH\in\dom[z,w]$ such that 
$\bV(z) = \bH\big(z,\bV(z)\big)$. 
If $\bE$ is open at $\big(\rho_\bT,\bT(\rho_\bT)\big)$
does it follow that 
$\bH$ is open at $\big(\rho_\bV,\bV(\rho_\bV)\big)$?

If so one would have an easy way of reducing the multi-singularity case 
of $\bT$ to the unique singularity case of $\bV$.
As mentioned in $\S\,$\ref{altern}, we were not able to prove this, 
but instead needed an additional hypothesis on $\bE$. 
Partly in order to avoid this extra hypothesis we used a detailed singularity 
analysis approach. 

\item[Q2]
Determine whether or not the $\Set_\bbM$ operators can be adjoined
to the standard operators used in this paper and still have the
universal law hold. (See $\S$\,\ref{Set sect}.)

A simple and interesting case to consider is that of identity
(0,1,...,m)-trees with generating function $\bT$ defined by 
$w = z + z\Set_{\{1,...,m\}}(w)$. Does $\bT$ satisfy the universal
law? Example 78 shows the answer is yes for $m=2$, but it seems 
the question is open for any $m \ge 3$.

\item[Q3]
Expand the theory to handle recursion equations $w = \bG(z,w)$ 
with $\bG$ having mixed sign coefficients.
\item[Q4]
Find large collections of classes 
 satisfying the universal law
(or any other law) that are recursively defined by 
systems of equations
\[
(**)\quad
\begin{array}{l l l}
w_1&=&\Theta_1\big(w_1,\ldots,w_k\big)\\
&\vdots&\\
w_k&=&\Theta_k\big(w_1,\ldots,w_k\big),
\end{array}
\]
where the $\Theta_i$ are multivariate operators. 

Classes defined by specifications using the standard operators 
that correspond to such a system of equations are called 
{\em constructible} classes by Flajolet and Sedgewick; 
the asymptotics of the case that the operators are {\em polynomial} has been
studied in \cite{Fl:Se}, Chapter VII, \emph{provided the dependency digraph 
has a single strong component}, 
and shown to satisfy the universal law $\pmb{(\star)}$.
\item[Q5]
The study of systems (**) is of particular interest to those investigating 
the behavior of monadic second order definable classes $\cT$ 
since every such class is a finite disjoint union of some of 
the $\cT_i$ defined by such a system (see Woods \cite{Wo}). 
In brief, Woods proved: \emph{every MSO class of trees with the same radius
as the whole class of trees satisfies the universal law.} However his results
seem to give very little for the MSO classes of smaller radius beyond the
fact that they have smaller radius.
Here is a plausible direction:

If $\cT$ is a class of trees defined by a MSO sentence, does it follow
that $\cT$ decomposes into finitely many $\cT_i$ such that each satisfies
a nice law on its support?

\item[Q6]
Among the MSO classes of trees perhaps the best known are the \emph{exclusion}
classes $\cT\,=\,\Excl(T_1,\ldots,T_n)$, defined by saying that certain trees 
$T_1,\ldots,T_n$ are not to appear as induced subtrees. 
The $T_i$ are called \emph{forbidden trees}\,.
A good example is
`trees of height $n$', defined by excluding a chain of height $n+1$\,;
or unary-binary trees defined by excluding the four-element tree of height 1.

Even restricting one's attention to the collection of classes 
$\cT\,=\,\Excl(T)$ defined by excluding a single tree 
offers considerable challenges to the development of a global theory of 
enumeration. 

Which of these classes are defined by recursion? 
Which of these obey the universal law? 
Given two trees
$T_1,T_2$\,, which of $\Excl(T_1),\Excl(T_2)$ has the greater radius (for its
generating series)? 

From Schwenk \cite{Schw} (1973) we know that if one excludes any limb from the
class of trees then the radius of the resulting class is larger than
what one started with.
Much later, in 1997,  Woods \cite{Wo} 
rediscovered a part of Schwenk's result in the context of logical
limit laws; this can be used to quickly show that the class of
free trees has a monadic second order 0--1 law. 
Aaron Tikuisis, an undergraduate at the University of Waterloo, has
determined which $\Excl(T)$ have radius $<1$. 
\item[Q7]
Find a method to determine the asymptotics of a class $\cU$ of free trees
obtained from a recursively defined class $\cT$ of rooted trees.
(See $\S\,$\ref{20 step}.)
\end{thlist}

\end{document}